\documentclass[11pt]{amsart}
\usepackage[T1]{fontenc} 
\usepackage{amsmath,amsthm,amsfonts,amssymb,bm,graphicx,xcolor}

\usepackage
{hyperref}
\hypersetup{colorlinks=true,citecolor=blue,linkcolor=blue,urlcolor=blue,
pdfstartview=FitH }

\theoremstyle{plain}
\newtheorem{theorem}{Theorem}
\newtheorem{cor}[theorem]{Corollary}
\newtheorem{prop}{Proposition}
\newtheorem{lemma}[prop]{Lemma}
\newtheorem{defi}{Definition}

\theoremstyle{definition}

\renewcommand{\geq}{\geqslant}
\renewcommand{\leq}{\le}

\newcommand{\Z}{\mathbb{Z}}
\newcommand{\Q}{\mathbb{Q}}
\newcommand{\N}{\mathbb{N}}
\newcommand{\R}{\mathbb{R}}
\newcommand{\C}{\mathbb{C}}
\newcommand{\F}{\mathbb{F}}

\newcommand{\A}{\mathbb{A}}

\newcommand{\bs}{\backslash}

\DeclareMathOperator{\Proj}{Proj}

\newcommand{\sgn}{\operatorname{sgn}}

\newcommand{\Stab}{\operatorname{Stab}}

\newcommand{\Vol}{\operatorname{Vol}}

\newcommand{\GL}{\operatorname{GL}}

\newcommand{\SL}{\operatorname{SL}}
\newcommand{\SO}{\operatorname{SO}}

\newcommand{\bdd}{\begin{center}\begin{tikzcd}}
\newcommand{\bd}{\begin{tikzcd}}
\newcommand{\edd}{\end{tikzcd}\end{center}}
\newcommand{\ed}{\end{tikzcd}}
\newcommand{\bdp}{\begin{center}\begin{tikzpicture}}
\newcommand{\edp}{\end{tikzpicture}\end{center}}
\newcommand{\bi}{\begin{itemize}}
\newcommand{\ei}{\end{itemize}}
\newcommand{\bt}{\begin{tikzpicture}}
\newcommand{\et}{\end{tikzpicture}}
\newcommand{\ba}{\[\begin{aligned}}
\newcommand{\ea}{\end{aligned}\]}
\newcommand{\bp}{\begin{pmatrix}}
\newcommand{\ep}{\end{pmatrix}}
\newcommand{\bv}{\begin{vmatrix}}
\newcommand{\ev}{\end{vmatrix}}
\newcommand{\bb}{\begin{bmatrix}}
\newcommand{\eb}{\end{bmatrix}}
\newcommand{\bB}{\begin{Bmatrix}}
\newcommand{\eB}{\end{Bmatrix}}
\newcommand{\bea}{\begin{enumerate}[label=(\alph*)]}
\newcommand{\ber}{\begin{enumerate}[label=(\roman*)]}
\newcommand{\ben}{\begin{enumerate}[label=(\arabic*)]}
\newcommand{\ee}{\end{enumerate}}

\newcommand{\mt}[4]{{\setlength{\arraycolsep}{0.5mm}\left(
\begin{smallmatrix}#1&#2\\#3&#4\end{smallmatrix}\right)}}

\numberwithin{equation}{section}
\def\I{\mathcal{I}}

\makeatletter
\def\Ddots{\mathinner{\mkern1mu\raise\p@
\vbox{\kern7\p@\hbox{.}}\mkern2mu
\raise4\p@\hbox{.}\mkern2mu\raise7\p@\hbox{.}\mkern1mu}}
\makeatother
\makeatletter
\DeclareRobustCommand\widecheck[1]{{\mathpalette\@widecheck{#1}}}
\def\@widecheck#1#2{%
    \setbox\z@\hbox{\m@th$#1#2$}%
    \setbox\tw@\hbox{\m@th$#1%
       \widehat{%
          \vrule\@width\z@\@height\ht\z@
          \vrule\@height\z@\@width\wd\z@}$}%
    \dp\tw@-\ht\z@
    \@tempdima\ht\z@ \advance\@tempdima2\ht\tw@ \divide\@tempdima\thr@@
    \setbox\tw@\hbox{%
       \raise\@tempdima\hbox{\scalebox{1}[-1]{\lower\@tempdima\box
\tw@}}}%
    {\ooalign{\box\tw@ \cr \box\z@}}}
\makeatother

\def\m{\mathsf{m}}

\begin{document}

\author{Valentin Blomer}
\address{Mathematisches Institut, Endenicher Allee 60, 53115 Bonn, Germany}
\email{blomer@math.uni-bonn.de}
 
\author{F\'elicien Comtat}
\email{comtat@math.uni-bonn.de}

 \title{Moments of symmetric square $L$-functions on ${\rm GL}(3)$}

\thanks{Both authors were supported by ERC Advanced Grant  101054336 and Germany's Excellence Strategy grant EXC-2047/1 -- 390685813. The first author was partially funded by the German Research Foundation -- Project-ID 491392403 -- TRR 358.}


\begin{abstract} We give an asymptotic formula with power saving error term for the twisted first moment of symmetric square $L$-functions on ${\rm GL}(3)$ in the level aspect. As applications, we obtain non-vanishing results as well as lower bounds of the expected order of magnitude for all even moments, supporting the random matrix model for a unitary ensemble. Besides the ${\rm GL}(3)$ Kuznetsov formula, the ingredients include detailed local computations at ramified places, including root numbers and orthonormalization of oldforms and Eisenstein series.
   \end{abstract}


\setcounter{tocdepth}{2}  \maketitle 

\maketitle

\section{Introduction}

\subsection{Twisted moments of prime level}

The analytic and statistical properties of the degree 6 family of symmetric square $L$-functions on ${\rm GL}(3)$ have received little attention so far, very much in contrast to the degree 6 family of ${\rm GL}(3) \times {\rm GL}(2)$ Rankin-Selberg $L$-functions. In this paper, we are interested in the central values  $L(1/2, \pi,  \text{sym}^2)$ for the family $\mathcal{B}(p)$ of cuspidal automorphic representations $\pi$ on ${\rm PGL}(3)$ over $\Bbb{Q}$ of odd prime level $p$. Classically this corresponds to  cuspidal newforms for the congruence subgroup $\Gamma_0(p) \subseteq {\rm SL}_3(\Bbb{Z})$ of index $p^2 + p + 1 \asymp p^2$. 



By Langlands' general conjectures it is expected that the symmetric square lift from ${\rm GL}(3)$ to ${\rm GL}(6)$ is functorial, but this is currently not known. In particular, the family of $L$-functions we consider in this paper has not yet been proved to be automorphic. As we will see, this makes certain parts of the analysis quite subtle.

Nevertheless, some analytic properties of the $L$-function $L(s, \pi, \text{sym}^2)$ are known. Its meromorphic continuation, functional equation and absolute convergence in $\Re s > 1$ follows from the factorization of the completed $L$-functions
\begin{equation}\label{factor}
\Lambda(s, \pi, \text{sym}^2) \Lambda(s, \tilde{\pi}) = \Lambda(s, \pi \otimes \pi)
\end{equation}
and the corresponding properties of Rankin-Selberg and standard $L$-functions. The hard part is to locate possible poles. Banks \cite[Theorem 7]{Ba} showed that $\Lambda^{(S)}(s, \pi, \text{sym}^2)$ is entire except for possible poles at $s = 0, 1$ where $S$ is the set of all places where $\pi$ ramifies as well as 2 and $\infty$. A pole at $s=1$ occurs if and only if $\pi$ is itself the symmetric square of a ${\rm GL}(2)$ representation (cf.\  \cite[p.\ 552]{PPS}).  Unfortunately we have   little control on the local factors at 2 and at $\infty$. At the current state of knowledge they could  have   poles anywhere in the strip $|\Re s| \leq 5/7$ by the best known bounds towards Ramanujan on ${\rm GL}(3)$. In particular there could be a pole at  $s=1/2$, in which case the problem of analyzing the central value $L(1/2, \pi,  \text{sym}^2)$ can become meaningless.  By a density theorem \cite[Theorem 4]{BBM}, this cannot happen too often, but the mere existence of such cases leads to major difficulties. We therefore slightly modify the $L$-function at the prime 2 by replacing the local factor $L_2(s, \pi, \text{sym}^2)$ with $L_2(s, \tilde{\pi}, \text{sym}^2)^{-1}$: we define
\begin{equation}\label{lstar}
L^{\ast}(s, \pi, \text{sym}^2) = \frac{L^{(2)}(s, \pi, \text{sym}^2)}{L_2(s, \tilde{\pi}, \text{sym}^2)}
\end{equation}
which \emph{vanishes} if $L_2(s, \tilde{\pi}, \text{sym}^2)$ has a pole at $s=1/2$. The possible poles of the local factor at infinity will be dealt with by choosing a suitable test function at the archimedean place, cf.\ Definition \ref{def1} below.

We will analyze this family of $L$-functions with the ${\rm GL}(3)$ Kuznetsov formula. As a relative trace formula, the Kuznetsov formula is based on the Fourier coefficients    $A_{\pi}(m, n)$ which are locally Schur polynomials in the Satake parameters of $\pi$ (in particular satisfying $A_{\pi}(1, 1) = 1$), cf.\ Subsection \ref{11}. A useful toolbox of independent interest for the spectral side of the Kuznetsov formula is provided in the subsequent subsections that investigate  Fourier coefficients at ramified primes, Fourier coefficients of oldforms and their various shifts, and root numbers expressed in terms of Fourier coefficients.  On the arithmetic side, we will have to analyze various kinds of exponential sums and   need strong bounds for Whittaker-type functions that appear as integrals transforms in the Kuznetsov formula. 

 Our main result is as follows.

\begin{theorem}\label{thm1}  There exist  constants $B, A_0, \delta > 0$  with the following property. Let $p$ be an odd prime and $\m_1, \m_2\in \Bbb{N}$ with $(\m_1\m_2, 2p) = 1$.  
Let $h \in \mathcal{H}(A_0)$ be a    test function as in Definition~\ref{def1} below. 
 Then there exists a constant $C$ depending only on $h$ such that  
$$\sum_{\pi \in  \mathcal{B}(p)} A_{\pi}(\m_1, \m_2) \frac{L^{\ast}(1/2, \pi, \text{{\rm sym}}^2)}{L(1, \pi, \text{{\rm Ad}})} h(\mu_{\pi}) = \delta_{\substack{\m_1 = \square\\ \m_2 = \square}} C \frac{p^2}{(\m_1\m^2_2)^{1/4}} + O_h(p^{2-\delta} (\m_1\m_2)^B).$$
\end{theorem}

A direct application of the asymptotic formula is a non-degeneracy result. 
\begin{cor}\label{cor2}  There exists a constant $c > 1$ with the following property. If $p > c$ is a prime, then not all numbers $L(1/2, \pi, \text{\rm sym}^2)$ for $\pi \in \mathcal{B}(p)$, $ \| \mu_{\pi}  \| \leq c$ are zero or undefined. 
\end{cor}
This follows trivially from the special case $\m_1 = \m_2 =1$ of Theorem \ref{thm1} and the choice of a test function $h$ which is positive on the tempered spectrum so that  $C > 0$, e.g.\ \eqref{testh}.  

\medskip

A more subtle application, which crucially needs the twist by $A_{\pi}(\m_1, \m_2)$ with polynomial dependence on $\m_1, \m_2$, are lower bounds for \emph{all} moments of the expected order of magnitude. This is based on a method of Rudnick and Soundararajan \cite{RS}. 

\begin{cor}\label{cor3}  Let $h$ be as in Theorem \ref{thm1} and non-negative on the automorphic spectrum, e.g.\ given by \eqref{testh}. Then for any fixed $k \in \Bbb{N}$ and $p$ prime we have
$$\sum_{\pi \in  \mathcal{B}(p)} \frac{|L^{\ast}(1/2,  \pi, \text{\rm sym}^2)|^{2k}}{L(1, \pi, \text{{\rm Ad}})} h(\mu_{\pi})  \gg_k p^2 (\log p)^{k^2}$$
if $p$ is sufficiently large in terms of $k$. 
\end{cor}

The exponent $k^2$ for the $2k$-th moment is in agreement with the random matrix model for a unitary ensemble of $L$-functions.\\


We now turn to a number of results that are needed as input for Theorem \ref{thm1} and may likely find applications elsewhere. The following subsections focus on the ramified Hecke algebra, symmetric square $L$-functions and orthogonalization of oldforms and Eisenstein series. 

\subsection{The   Hecke algebra}\label{11}

We recall that $A_{\pi}(m, n)$ is multiplicative in both arguments and on prime powers given by
\begin{equation}\label{schur}
A_{\pi}(p^k, p^{\ell}) = \det \left(\begin{smallmatrix} \alpha_p^{k+\ell + 2} & \beta_p^{k+\ell + 2} & \gamma_p^{k + \ell + 2}\\ \alpha_p^{k+1} & \beta_p^{k+1} & \gamma_p^{k+1}\\ 1 & 1 & 1\end{smallmatrix}\right)/\det \left(\begin{smallmatrix} \alpha_p^{2} & \beta_p^{  2} & \gamma_p^{ 2}\\ \alpha_p  & \beta_p  & \gamma_p \\ 1 & 1 & 1\end{smallmatrix}\right)
\end{equation}
 in terms of the Satake parameters $\alpha_p, \beta_p, \gamma_p$ of $\pi$ at $p$. If $p$ is unramified, then
\begin{equation}\label{proper}
\begin{split}
 & \alpha_p\beta_p \gamma_p = 1 \quad  \text{(trivial central character)},\\
  & \{\alpha^{-1}_p, \beta_p^{-1}, \gamma^{-1}_p\} = \{\overline{\alpha_p}, \overline{\beta_p}, \overline{\gamma_p}\} \quad \text{(unitarity)},\\
  &  |\alpha_p|, |\beta_p|, |\gamma_p| \leq p^{\theta_3} \quad \text{(bounds towards Ramanujan)}
  \end{split}
\end{equation} 
with
$$\theta_3 = 5/14$$
by the Kim-Sarnak bound, whose proof is summarized in \cite{BB}. 
The formula \eqref{schur} together with multiplicativity implies the usual Hecke relations. We mention in particular (\cite[Lemma 2.1]{BL}, \cite[Theorem 6.4.11]{Go})
\begin{equation}\label{hecke}
A_{\pi}(m, n_1 n_2 ) = \sum_{\substack{abc = n_1\\ b \mid mc, c \mid n_2}} \mu(b) \mu(c) A_{\pi}\Big( \frac{mc}{b}, \frac{n_2}{c}\Big) A_{\pi}(1, a)
\end{equation}
and
\begin{equation}\label{hecke1}
  A_{\pi}(n, 1) A_{\pi}(m_1, m_2) = \sum_{\substack{d_0d_1d_2 = n\\ d_1 \mid m_1, d_2 \mid m_2}} A\Big(\frac{m_1d_0}{d_1}, \frac{m_2 d_1}{d_2}\Big).
\end{equation}
The second formula in \eqref{proper} implies for instance $A_{\pi}(m, n) = \overline{A_{\pi}(n, m)}$ as long as $m, n$ have only unramified prime divisors. We will frequently use the bound
$$A_\pi(m, n) \ll (mn)^{\theta_3 + \varepsilon}.$$

If $p$ is ramified, we still define $A_{\pi}(p^k, p^{\ell})$ by \eqref{schur}, but \eqref{proper} no longer holds.  

\begin{lemma}\label{ram} Let $\pi$ be a unitary automorphic representation for ${\rm PGL}(3)$ of prime level $p$. Then there exists an imaginary number $\rho\in i\Bbb{R}$ such that $$
\alpha_p = p^{-1/2 - \rho}, \quad \beta_p = p^{2\rho}, \quad \gamma_p = 0.$$
In particular we have the Hecke relations
$A_{\pi}(p^j, 1) = p^{j(-1/2 - \rho)}$,  $A_{\pi}(p^j, p^i) =  A_{\pi}(p^j,1)A_{\pi}(1, p^i)$ for $i, j \geq 0$, and the last bound in \eqref{proper} remains valid, even with $\theta_3 = 0$. 
\end{lemma}

 \subsection{Symmetric square $L$-functions}\label{12}
 
 Next we consider the Dirichlet series of symmetric square $L$-functions on ${\rm GL}(3)$ of prime (or squarefree) level, their conductor and root number as well as regularity properties. 
    
 \begin{prop}\label{prop1}
Let $N \in \Bbb{N}$ be squarefree and $\pi$ a cuspidal representation of level $N$ and spectral parameter $\mu$. 
For $\Re s > 1$ we have 
 $$L(s,\pi,  \text{{\rm sym}}^2) = \sum_{\substack{a, b, c\\ (a, N) = 1}}  \frac{A_{\pi}(b^2, c^2)}{(a^3b^2c)^s}.$$
 In particular, the local $L$-factor at a prime $p \mid N$ is given by
 $$L_p(s, \pi,  \text{\rm sym}^2) = \Big(1 - \frac{p^{-2\rho}}{p^{1 + s}}\Big)^{-1}\Big(1 - \frac{p^{-\rho}}{p^{1/2 + s}}\Big)^{-1}\Big(1 - \frac{p^{4\rho}}{p^{s}}\Big)^{-1}$$
 with $\rho \in i\Bbb{R}$ as in Lemma \ref{ram}. The local factor at infinity equals
 $$L_{\infty}(s, \pi, \text{\rm sym}^2) =   \prod_{1 \leq i \leq j \leq 3}  \Gamma_{\Bbb{R}} (1/2 + u - \mu_{  i} - \mu_{  j}), \quad  \Gamma_{\Bbb{R}}(s) = \Gamma(s/2)\pi^{-s/2}.$$
  \end{prop}

 \begin{prop}\label{prop2} Let $N \in \Bbb{N}$ be squarefree and $\pi$ a cuspidal representation of level $N$. 
 The  arithmetic conductor of $L(s, \pi, \text{\rm sym}^2)$ is $N^3$ and 
 the root number $\epsilon = \epsilon(\pi, \text{\rm sym}^2, 1/2)  
 $ is given by
 $$\epsilon(\pi, \text{\rm sym}^2, 1/2)  = \prod_{p \mid N} \big(-p^{3/2} \overline{A_{\pi}(p^3, 1)}\big) = \prod_{p \mid N} \left( - p^{1/2}\overline{A_{\pi}(p, p)} + \frac{1}{p^{1/2}}\right).$$
 \end{prop}

  A direct consequence of Proposition \ref{prop1} is the following
   \begin{lemma}\label{polefree} Let $N \in \Bbb{N}$ be squarefree and $\pi$ a cuspidal representation of level $N$. 
  Then the partial $L$-function $L^{(2)}(s, \pi, \text{\rm sym}^2)$ has (at most simple) poles at most at $s=0, 1$.
  \end{lemma}
  
  \textbf{Proof.} We know from \cite[Theorem 7]{Ba} that $L^{(2N)}(s, \pi, \text{\rm sym}^2)s(1-s)$  is entire. The functional equation implies
\begin{displaymath}
\begin{split}
&L^{(2N)}(s, \pi, \text{\rm sym}^2)s(1-s)N^{3s/2} \prod_{v \mid 2N \infty} L_v(s, \pi,  \text{\rm sym}^2)  \\
& = \epsilon L^{(2N)}(1-s, \pi, \text{\rm sym}^2)s(1-s) N^{3(1-s)/2}\prod_{v \mid 2N \infty} L_v(1-s, \pi,  \text{\rm sym}^2) .
\end{split}
\end{displaymath}
For $p \mid N$, by Lemma \ref{ram} the   local $L$-factor $ L_p(s, \pi,  \text{\rm sym}^2)$   can have poles at most in $\Re s \leq 0$ 
where the right hand side is holomorphic.
 
 \subsection{Spectral decomposition and orthonomalization}\label{13}
 
Our next results concern the image of  ``oldforms''  of level 1 in $L^2(\Gamma_0(p)\bs {\rm SL}_3(\Bbb{R})/{\rm SO}(3))$, in particular with respect to an orthonormal basis of the  translates of oldforms.  We equip all spaces $L^2(\Gamma_0(N)\bs {\rm SL}_3(\Bbb{R})/{\rm SO}(3))$, $N \in \Bbb{N}$, with a probability Haar measure so that the embedding $$L^2({\rm SL}_3(\Bbb{Z})\bs {\rm SL}_3(\Bbb{R})/{\rm SO}(3))\hookrightarrow L^2(\Gamma_0(p)\bs {\rm SL}_3(\Bbb{R})/{\rm SO}(3))$$ is an isometry.

 \begin{prop}\label{prop3} Let $p$ be a prime. Let $\pi$ be a representation of level 1, and let $f$ in $\pi$ be a spherical vector.  
 The space of oldforms generated by $f$ for the group $\Gamma_0(p)$ is generated by three forms $T_0f$, $T_1f$, $T_2f$ with Fourier coefficients
\begin{equation}\label{oldFour}
\begin{split}
&A^{(0)}_{\pi}(m, n) = A_{\pi}(m, n), \\  
& A^{(1)}_{\pi}(m, n) = \sqrt{p}\Big(A_{\pi}\Big(\frac{m}{p}, pn\Big) + A_{\pi}\Big(m, \frac{n}{p}\Big)\Big), \\
& A^{(2)}_{\pi}(m, n) = pA_{\pi}\Big(\frac{m}{p}, n\Big)
 \end{split}
 \end{equation}
respectively,  with the convention $A_{\pi}(m, n) = 0$ if $m$ or $n$ are not integral. The corresponding Gram matrix of pairwise inner products is given by
 $$\left(\begin{matrix} 1 &\frac{p^{5/2} - p^{1/2} }{p^3 - 1}A_{\pi}(1, p) & \frac{p^2  - p }{p^3 - 1}  A_{\pi}(p, 1)\\ \frac{p^{5/2} - p^{1/2} }{p^3 - 1}A_{\pi}(p, 1)   & 1+\frac{ p^2 - p}{p^3 - 1}|A_{\pi}(1, p)|^2&\frac{p^{5/2} - p^{1/2} }{p^3 - 1}A_{\pi}(1, p)  \\ \frac{p^2  - p }{p^3 - 1}  A_{\pi}(1, p)& \frac{p^{5/2} - p^{1/2} }{p^3 - 1}A_{\pi}(p, 1)& 1 \end{matrix}\right) \| f \|^2.$$
 \end{prop}

\textbf{Remark:}  It might be a little counter-intuitive that the Fourier coefficients of $T_2f$ are of size $p$, but only every $p$-th Fourier coefficient is picked up. Square-root cancellation in the Fourier expansion would suggest to put a weight $p^{1/2}$ on the Fourier coefficients as in the case of $T_1f$. The reason is that the ${\rm GL}(3)$ Fourier expansion is more   complicated and contains an additional ${\rm SL}(2)$-sum. The above normalization is consistent with square-root cancellation in the entire Fourier expansion (which puts more weight on the second entry of the Fourier expansion). 

\begin{cor}\label{cor4} With the assumptions and notation of Proposition \ref{prop3}, the space of oldforms generated by $f$ for $\Gamma_0(p)$ has an orthonormal basis $S_0f, S_1f, S_2f$ whose Fourier coefficients are given by $B_{\pi}^{(j)}(m, n) \| f \|^{-1}$, $j = 0, 1, 2$, respectively, where
\begin{equation}\label{orth}
\begin{split}
B_{\pi}^{(0)}(m, n) = &A_{\pi}^{(0)}(m, n), \\  
B_{\pi}^{(1)}(m, n) =& (1 + c_{11} )A^{(1)}_{\pi}(m, n)  + c_{10} A^{(0)}_{\pi}(m, n),\\
B_{\pi}^{(2)}(m, n) =& (1 + c_{22}) A^{(2)}_{\pi}(m, n) + c_{21}A^{(1)}_{\pi}(m, n) + c_{20} A^{(0)}_{\pi}(m, n)
 \end{split}
 \end{equation}
for suitable constants $c_{ij} \ll p^{\theta_3 - 1/2}$ depending only $A_{\pi}(1, p)$ and $A_{\pi}(p, 1)$. 
\end{cor}

\textbf{Proof.} 
The entries of the Gram matrix in Proposition \ref{prop3} are 
$$\left(\begin{matrix} 1 & O(p^{\theta_3 - 1/2 })  & O(p^{\theta_3 - 1 }) \\O(p^{\theta_3 - 1/2  })   & 1 + O(p^{2\theta_3 - 1  }) &O(p^{\theta_3 - 1/2  })  \\  O(p^{\theta_3 - 1  })& O(p^{\theta_3 - 1/2  }) & 1 \end{matrix}\right) \| f \|^2 = {\rm id} \cdot \| f\| ^2 + O(\| f \|^2p^{\theta_3 - 1/2 }).$$
We can now run the standard Gram-Schmidt algorithm to orthonormalize $\{T_0f, T_1f, T_2f\}$.\\ 

Next we  discuss the spectral decomposition of the space of unitary Eisenstein series of level $p$. We start with the minimal Eisenstein series. They come from ``oldforms'' of level one, and we recall that such Eisenstein series are parametrized by $s = (s_1, s_2) \in (i\Bbb{R})^2$ with local Satake parameters
$$\alpha_p = p^{s_1}, \quad \beta_p = p^{s_2}, \quad \gamma_p = p^{-s_1-s_2}$$
satisfying \eqref{proper}. We define $A_{(s_1, s_2)}(m, n)$ by \eqref{schur} and multiplicativity. 
As in the cuspidal case each oldform generates a 3-dimensional subspace, and the relation between Fourier coefficients and Gram matrices turns out to be identical. Thus for $j=0, 1, 2$ we define $A^{(j)}_{(s_1, s_2)}(m, n)$ in terms of $A_{(s_1, s_2)}(m, n)$ as in \eqref{oldFour} and $B_{(s_1, s_2)}^{(j)}(m, n)$ in terms of $A^{(j)}_{(s_1, s_2)}(m, n)$ as in \eqref{orth}. 
\begin{prop}\label{prop4}
With the notation as above,  the space of minimal unitary Eisenstein series is parametrized by $(j, (s_1, s_2)) \in \{0, 1, 2\} \times (i\Bbb{R})^2$ and the corresponding Eisenstein series has  Fourier coefficients $B_{(s_1, s_2)}^{(j)}(m, n)$ where the constants $c_{ij}$ satisfy the bound $c_{ij} \ll p^{  -1/2}$. 
\end{prop} 

We now turn to (generic) maximal Eisenstein series of level $p$. These are parametrized by a Maa{\ss} form $\phi$ for ${\rm PGL}(2)$ with Hecke eigenvalues $\lambda_{\phi}(n)$ and spectral parameter $t_{\phi}$, say, and an  imaginary parameter $s\in i\Bbb{R}$. If $p$ is a prime and  $a_p, b_p$ are the local Satake parameters of $\phi$ (satisfying $|a_p|, |b_p| \leq p^{\theta_2}$ with $\theta_2 = 7/64$), then we define the local Satake parameters
\begin{equation}\label{satmax}
\alpha_p = a_pp^{s} \quad \beta_p = b_p p^{s}, \quad \gamma_p = p^{-2s}.
\end{equation}
If $\phi$ is unramified at $p$, they satisfy \eqref{proper}. If $\phi$ is ramified at $p$, then $\{a_p, b_p\}  = \{p^{-1/2}, 0\}$, and the numbers in \eqref{satmax} are, up to permutation, of the form as in Lemma \ref{ram}. Correspondingly we define $A_{(\phi, s)}(m, n)$ by \eqref{schur} and multiplicativity. If $\phi$ has level 1, we define as above $A^{(j)}_{(\phi, s)}(m, n)$  and $B_{(\phi, s)}^{(j)}(m, n)$ using \eqref{oldFour} and  \eqref{orth}. Running now the same argument, we obtain




\begin{prop}\label{prop5}
With the notation as above,  the space of maximal generic unitary Eisenstein series is parametrized by

a) pairs $(\phi, s)$ of cuspidal Maa{\ss} newforms of level $p$ with Fourier coefficients $A_{(\phi, s)}(m, n)$. 

b) triplets  $(j, \phi, s)$ with $j \in \{0, 1, 2\}$ and $\phi$ of level 1,  with corresponding  Fourier coefficients $B_{(\phi, s)}^{(j)}(m, n)$ as in \eqref{orth} with $c_{ij} \ll p^{\theta_2 - 1/2}$. 
\end{prop}


 \section{Local computations I: symmetric square $L$-functions}
 
 In this section we prove Lemma \ref{ram} and the propositions in  Subsection  \ref{12}.\\

 We denote by $\nu$ the standard absolute value on $\Q_p$. Given a smooth irreducible representation~$\pi_p$ of~${\rm GL}_n(\Q_p)$,
there exists a non-negative integer $a$ such that $\pi$ has 
a non-trivial vector fixed by
$$
\Gamma_0^{(n)}(p^a):=\left(\begin{matrix} 
M_{n-1,n-1}(\Z_p) & M_{n-1,1}(\Z_p)\\ 
p^aM_{1,n-1}(\Z_p) & 1+p^a\Z_p \end{matrix}\right)
\cap {\rm GL}_n(\Q_p).$$
By definition the least such integer $a=a(\pi_p)$ is the conductor of $\pi_p$.\\

  \textbf{Proof} of~Lemma \ref{ram}. Consider the local component $\pi_p$ of $\pi$, viewed as a representation of $\GL_3(\Q_p)$ with trivial central character. By assumption, we have $a(\pi_p)=1$.  Since supercuspidal representations of $\GL_3$ have conductor at least 3, $\pi$ cannot be supercuspidal.  Hence by the Bernstein-Zelevinsky classification, $\pi$ is the unique irreducible quotient of a representation of one of the following forms:
    \begin{enumerate}
        \item $\chi \times \chi \nu \times \chi \nu^2$ for some character $\chi$,
        \item $\chi_1 \times\chi_2\rm{St}_2$, where $\chi_1 \neq \chi_2\nu^{-3/2},\chi_2\nu^{3/2}$,
        \item  $\chi \times \sigma$ where $\sigma$ is a supercuspidal representation of $\GL_2$,
        \item $\chi_1 \times \chi_2 \times \chi_3$, where $\chi_i\chi_j^{-1} \neq \nu^{\pm1}$ for $i \neq j$.
    \end{enumerate}
      First, we show that only case $(2)$ is possible. Indeed, the irreducible quotient of $(1)$ is a twist of ${\rm St}_3$, which has conductor $2$.
    Next, the representation $(3)$ is irreducible and hence has conductor $a(\sigma)+a(\chi)$, but $a(\sigma) \ge 2$. Finally, the representation $(4)$ is also irreducible.
    Now if $\pi_p=\chi_1 \times \chi_2 \times \chi_3$ 
    then we have $1=a(\pi_p)=a(\chi_1)+a(\chi_2)+a(\chi_3)$ hence exactly one
    of the three characters is ramified. Thus $\chi_1\chi_2\chi_3$ is ramified.
    But  $\chi_1\chi_2\chi_3=1$ since the central character is trivial, which is a contradiction.

      Now, $\chi_1 \times\chi_2\rm{St}_2$ is irreducible, and thus it has conductor $a(\chi_1)+a(\chi_2\rm{St}_2)$.
    This forces $a(\chi_1)=a(\chi_2)=0$. Comparing the central characters, we obtain $\chi_1=\chi_2^{-2}$.   

    Finally, since $\pi_p$ is unitary, the contragradiant $\tilde \pi_p$ is isomorphic to the complex conjugate $\overline{\pi}_p$ of $\pi_p$.
But 
$$\tilde \pi_p=\chi_2^{-1}\rm{St}_2 \times \chi_2^2 \quad \text{and}\quad \overline{\pi}_p=\overline{\chi}_2\rm{St}_2 \times \overline{\chi}_2^{-2}.$$
By the Bernstein-Zelevinsky classification, these two representations are isomorphic if and only if $\chi_2^{-1}=\overline{\chi}_2$. Since moreover $a(\chi_2)=0$, 
there exists an imaginary number $\rho \in i\R$ such that
$\chi_2=\nu^\rho$. 

To show that $\alpha_p=p^{-1/2-\rho},$ $\beta_p=p^{2\rho}$, $\gamma_p=0$ are the Satake parameters of $\pi_p$ as defined in~\S\ref{11}, a simple calculation shows that 
$$L_p(s,\pi)=(1-\alpha_p p^{-s})^{-1}(1-\beta_p p^{-s})^{-1},$$
and the result now follows from \cite[Theorem 4.1]{Miy}.\\

 \textbf{Proof} of Proposition \ref{prop1}.  This is a purely local computation by comparing Euler factors on both sides. Let first $p \nmid N$ be an unramified prime, and let $\alpha_p, \beta_p, \gamma_p$ denote the Satake parameters of $\pi$ at $p$.  We express $A_{\pi}(p^k, p^{\ell})$ 
 in terms of Schur polynomials by \eqref{schur} and compute
 $$\sum_{i, j, k} \frac{A_{\pi}(p^{2j}, p^{2i})}{p^{(i+2j+3k)s}}$$
 by first evaluating the determinant in the numerator and then computing a finite number of  geometric series. A straightforward computation gives
 $$\frac{1 - (\alpha_p\beta_p\gamma_p)^2 p^{-3s}}{1 - p^{-3s}} \Big( \Big(1 - \frac{\alpha_p^2}{p^s}\Big) \Big(1 - \frac{\beta_p^2}{p^s}\Big) \Big(1 - \frac{\gamma_p^2}{p^s}\Big)\Big(1 - \frac{\alpha_p\beta_p}{p^s}\Big) \Big(1 - \frac{\alpha_p\gamma_p}{p^s}\Big) \Big(1 - \frac{\beta_p\gamma_p}{p^s}\Big)\Big)^{-1}.$$
Since $\alpha_p\beta_p\gamma_p = 1$, this coincides with the $p$-th Euler factor of $L(s, \text{{\rm sym}}^2, \pi) $.

If $p \mid N$, we apply the same argument without the $i$-sum, which removes the factor $1 - p^{-3s}$ in the denominator of the previous display, and we recall from Lemma \ref{ram} that $\alpha_p\beta_p\gamma_p = 0$. Again we obtain the $p$-th Euler factor of $L(s, \text{{\rm sym}}^2, \pi) $.

The local factor at infinity follows from the factorisation~(\ref{factor}) and the corresponding known local factors for $L(s, \pi \otimes \pi)$ and $L(s, \tilde{\pi})$.\\

 \textbf{Proof} of Proposition \ref{prop2}. The factorisation~(\ref{factor})
 implies the same factorisation for the sign of the functional equations. Thus 
 $$\epsilon(\pi, \text{\rm sym}^2, 1/2)=\prod_{p \mid N}
 \frac{\epsilon_p(\pi \times \pi, 1/2)}{\epsilon_p(\tilde{\pi}, 1/2)}.$$
  To compute the local epsilon factors on the right hand side of previous display, one may for instance use the local Langlands correspondence. 
   Let $p \mid N$.
 By the proof of Lemma~\ref{ram}, we know that
 $\pi_p=\nu^{-2\rho} \times \nu^{\rho} \rm{St}_2$.
 So its image by the local Langlands correspondence is
 given by $\rho_\pi=|\cdot|^{-2\rho} \oplus |\cdot|^{\rho-\frac12} \otimes \rm{Sp}_2$,
 where $\rm{Sp}_2$ is the special representation of dimension~$2$.
 Explicitly, $\rho_\pi=(r,N)$ is a representation of dimension~$3$
 that can be realised on a vector space $V=\langle e_0, e_1, e_2 \rangle$ by
     \begin{align*}
          &r(w)(e_0)=|w|^{-2\rho}e_0, &N(e_0)=0,\\
        &r(w)(e_1)=|w|^{\rho-\frac12}e_1, &N(e_1)=e_2,\\
    &r(w)(e_2)=|w|^{\rho+\frac12}e_2, &N(e_2)=0. 
    \end{align*} 
    for all $w \in W_{\Q_p}$.
 By definition of the local Langlands correspondence we have
 $$\epsilon_p(\pi \times \pi,s)=\epsilon(\rho_\pi \otimes \rho_\pi, s).$$
 In turn, the representation $\rho_\pi \otimes \rho_\pi$ is realised on the 9-dimensional vector space $W=V \otimes V$
 by $$(r \otimes r)(w) (e_i \otimes e_j) = (r(w) (e_i))\otimes (r(w) (e_j))=|w|^{a_{ij}} e_i \otimes e_j,$$
 $$(N \otimes N)(e_i \otimes e_j)=
 N (e_i) \otimes e_j + e_i \otimes N(e_j)= \begin{cases}
     e_2 \otimes e_1 + e_1 \otimes e_2 & \text{ if } i=j=1,\\
     e_2 \otimes e_j & \text{ if } i=1\neq j,\\
     e_i \otimes e_2 & \text{ if } i \neq 1=j,\\
     0 & \text{ otherwise.}
 \end{cases}$$
 By definition, we have 
$$
\epsilon(\rho_\pi \otimes \rho_\pi,s)=\epsilon(|\cdot|^{s}(r \otimes r))\det(-|\cdot|^{s}\Phi|_{W^I/W_{N \otimes N}^I}),
$$
 where
 $\Phi$ is a geometric Frobenius element, 
 $W^I$ is the subspace of $r$
 that is invariant by the inertia subgroup $I$, and $W_{N \otimes N}^I$ is the kernel
 of $N \otimes N$ in $W^I$.
  By multiplicativity in exact sequences, we have 
  $$\epsilon(|\cdot|^{s}(r \otimes r))=\prod_{i,j}\epsilon(|\cdot|^{s+a_{ij}})=1.$$
Therefore 
 $$\epsilon_p(\pi \times \pi, s)=\det(-|\cdot|^{s}\Phi|_{W^I/W_{N \otimes N}^I})=p^{-2\rho-4(s-1/2)}.$$
 Similarly, one may compute
 $$\epsilon_p(\tilde{\pi}, s)=\epsilon(\rho_{\tilde{\pi}_p}, s)=-p^{\rho-(s-1/2)}$$
for $p \mid N$, and thus 
 $$\frac{\epsilon_p(\pi \times \pi, 1/2)}{\epsilon_p(\tilde{\pi}, 1/2)}=-p^{-3\rho},$$
 which matches the factor on the right hand side in Proposition~\ref{prop2}, as one can check using Lemma~\ref{ram} and~\eqref{schur}. 

 \section{Eisenstein series}

 In this section we explicate the parametrisation of Eisenstein series. It is best to work adelically.  Let $\A$ be the ring of ad\`eles of $\Q$. For each place $v$ we define a 
 compact subgroup of $\GL_3(\Q_v)$ by setting $$K_p=\GL_3(\Z_p) \supseteq \Gamma_0^{(3)}(p) = \left\{\left(\begin{smallmatrix}
\Bbb{Z}_p & \Bbb{Z}_p & \Bbb{Z}_p \\
 \Bbb{Z}_p & \Bbb{Z}_p  & \Bbb{Z}_p \\
p \Bbb{Z}_p & p\Bbb{Z}_p  &\Bbb{Z}_p
\end{smallmatrix}\right)\right\} \cap K_p$$
if $p$ is a prime and $K_\infty=\SO_3(\R)$. 
We also define $\Gamma_0(p)=K_\infty \prod_p \Gamma_0^{(3)}(p)$.
This is the adelic counterpart of the classical Hecke congruence subgroup (which is denoted by the same symbol in the rest of the text).

We then set $K=\prod_{v}K_v$.
For $(s_1,s_2) \in \C^2$ define
$$I_{(s_1,s_2)}\left(\left(\begin{smallmatrix}
y_1 &  * &   * \\
  & y_2  & * \\
  &   & y_3
\end{smallmatrix}\right)\right)=|y_1|^{s_1}|y_2|^{s_2}|y_3|^{-s_1-s_2}$$ 
for all $y_1,y_2,y_3 \in \A^\times,$ and
$I_{(s_1,s_2)}(bk)=I_{(s_1,s_2)}(b)$ if $b \in \left(\begin{smallmatrix}
* &  * &   * \\
  & * & * \\
  &   & *
\end{smallmatrix}\right) \cap \GL_3$ and $k \in K$.

Let $P$ be a standard parabolic subgroup of  $\GL_3$ and write its Levi decomposition $P=NM$. 
Let $A$ be the centre of $M$.
Let $\sigma$ be a representation occurring in the discrete spectrum of 
$L^2(M(\Q)A^+(\R) \bs M(\A))$.
Langlands' theory attaches Eisenstein series to certain representations $\mathcal{I}_P(\sigma_s)$  induced from $\sigma$.
The representation $\mathcal{I}_P(\sigma_s)$ acts on the space $\mathcal{H}_P(\sigma)$ that consists of functions  $\phi: N(\A)M(\Q)A^+(\R) \bs \GL_3(\A) \to \C$
such that 
\begin{itemize}
    \item for all $x \in \GL_3(\A)$ the function
$\phi_x: m \to \phi(mx)$ belongs to $\sigma$, and
\item $\|\phi\|^2:=\int_K \int_{M(\Q)A^+(\R) \bs M(\A)} |\phi(mk)|^2 \, dm \, dk < \infty.$
\end{itemize} Moreover if $P$ and $P'$ are two parabolic subgroups that are associated to each other, then
their Eisenstein series define the same functions on $\rm{GL}_3(\Q) \bs \rm{GL}_3(\A)$.  More precisely, there are intertwining operators $M(w,s)$ between the representations induced from $P$ and $P'$, and these intertwining operators are unitary when $s$ is unitary, and they satisfy  $E(x,\phi,s)=E(x,M(w,s)\phi,s)$. 
Thus we only need to consider one representative in each association class of parabolic subgroups.

If $P=P_{111}$ is the minimal parabolic and $s=(s_1,s_2) \in \C^2$, we let $\mathbf{s}_P=(1+s_1,s_2)$.
If $P=P_{21}$ and $s \in \C$, we let $\mathbf{s}_P=(\frac12+s,\frac12+s)$.
    Note that if $g \in P_{21}$ then $$I_{\mathbf{s}_{P_{21}}}(g)=|\det(\mathfrak{m}_2(g))|^{\frac12+s}|\mathfrak{m}_1(g)|^{-1-2s},$$
where $\mathfrak{m}_1$ and $\mathfrak{m}_2$ are the projection on $\GL_1$ and $\GL_2$ respectively in the Levi decomposition
$P=NM$ and $M\simeq \GL_1 \times \GL_2$.

With these notations, the action of $\mathcal{I}_P(\sigma_s)$ on $\phi \in \mathcal{H}_P(\sigma)$ is given in both cases by
\begin{equation}\label{Ip}
\mathcal{I}_P(\sigma_s,y)\phi(x)=\phi(xy)I_{\mathbf{s}_P}(xy)
I_{-\mathbf{s}_P}(x),
\end{equation}
and the associated Eisenstein series satisfy
\begin{equation}\label{Intertwin}
    E(xy,\phi,s)=E(x,\mathcal{I}_P(\sigma_s,y)\phi,s).
\end{equation}

\begin{lemma}\label{Isom}
Fix a set of representatives $X$ of the double quotient 
$(P \cap K) \bs K /\Gamma_0(p)$.
For each $k \in X$, define
$M_k=\Proj_P^M(P \cap k \Gamma_0(p)k^{-1})$. 
Then there is an isomorphism
\begin{equation}\label{isomorphism}
\begin{split}
\I_P(\sigma_s)^{\Gamma_0(p)} & \to \bigoplus_{k \in X} \sigma^{M_k}\\
\phi & \mapsto (\phi_k)_{k \in X}.
\end{split}
\end{equation}
\end{lemma}
\textbf{Proof.}
   From equations~(\ref{Intertwin}) and~(\ref{Ip}),
   the Eisenstein series $E(\cdot, \phi,s)$ is 
   $\Gamma_0(p)$-invariant if and only if $\phi$ is $\Gamma_0(p)$-invariant.
By the Iwasawa decomposition $G(\A)=PK$, such a function $\phi$ is determined by the family $(\phi_k)_{k \in K/\Gamma_0(p)}$, which must in addition satisfy$$\phi(mgk)=\phi(\underbrace{(mn_1m^{-1})}_{\in N(\A)}mm_1k)\\
=\phi(mm_1k),$$
 for all $g=n_1m_1 \in P \cap K$ and for all $m \in M(\A)$, and hence
\begin{equation}\label{inducingdata}
\phi_{gk}=\sigma(m_1)\phi_k.
\end{equation}
Thus $\phi$ is completely determined by $(\phi_k)_{k \in X}$.
Moreover, by~(\ref{inducingdata}) the function $\phi_k$ is an $M_k$-invariant vector in the space of $\sigma$ for each $k \in X$, where
\begin{align*}
    M_k&=\Proj_P^M(\Stab_{P \cap K}(k \Gamma_0(p))=\Proj_P^M(P \cap k \Gamma_0(p)k^{-1}).
\end{align*}
Conversely, if each $\phi_k$ belongs to $\sigma^{M_k}$, then
we can define a $\Gamma_0(p)$-invariant function on
$N(\A)M(\Q)A^+(\R) \bs \GL_3(\A)$ 
by setting 
\begin{equation}\label{Inducingdata2}
    \phi(nmk\gamma)=\phi_k(m)
\end{equation} for all $n \in N(\A)$, $m \in M(\A)$ and for all $\gamma \in \Gamma_0(p)$. \\

Our aim is to determine an orthonormal basis of $\I_P(\sigma_s)^{\Gamma_0(p)}$.
Note that by~(\ref{Inducingdata2}) we have 
\begin{equation}\label{InnerProduct}
    \langle \phi, \psi \rangle =\int_K \langle \phi_k, \psi_k\rangle_M \, dk = \sum_{k \in X}v_k \langle \phi_k, \psi_k\rangle_M
\end{equation}
for $\phi, \psi \in \I_P(\sigma_s)^{\Gamma_0(p)}$ 
where $v_k=\Vol_K((P \cap K) k \Gamma_0(p) )$ and
$$\langle \phi_k, \psi_k\rangle_M= \int_{M(\Q)A^+(\R) \bs M(\A)} \phi_k(m) \overline{\psi_k}(m) \, dm.$$
At all places $v \neq p$, the double quotient $(P \cap K_v) \bs K_v / \Gamma_0(p)$ is trivial and thus the local component of $M_k$ is given by
$M \cap K_v$. We now turn to the calculation of the double quotient 
$(P \cap K_p) \bs K_p / \Gamma_0(p)$, the 
local component $M_{k,p}$ of $M_k$ at $p$, and the corresponding volumes~$v_{k,p}$.

\begin{lemma}\label{reps1}
Let $P=P_{111}$.
The set $$X_p=\left\{1,
w_1=\left(\begin{smallmatrix}
1 &   &    \\
  &   &  1\\
  & 1  & 
\end{smallmatrix}\right),w_2=
\left(\begin{smallmatrix}
 &   &  1  \\
 1 &   &  \\
  &  1 & 
\end{smallmatrix}\right)
\right\}$$ is a set of representatives of
$P(\Z_p) \bs \GL_3(\Z_p) / \Gamma_0(p) $.
The corresponding $M_{k,p}$ are given by $M_k=M \cap K$,
and the volumes are given by  $v_{1,p}=\frac{p-1}{p^3-1}$,
    $v_{w_1,p}=pv_1$ and $v_{w_2,p}=p^2v_1$.
\end{lemma}
\textbf{Proof.}
    We have 
    $P(\Z_p) \bs \GL_3(\Z_p) / \Gamma_0(p)= P(\F_p) \bs \GL_3(\F_p) / P_{21}(\F_p).$
    By the Bruhat decomposition over $\F_p$, a set of representatives 
    of
    $$P(\F_p) \bs \GL_3(\F_p) / P_{21}(\F_p)$$
    is given by the Weyl group, modulo (on the right) the subgroup generated by
    $\left(\begin{smallmatrix}
 &  1 &    \\
 1 &   &  \\
  &   & 1
\end{smallmatrix}\right),$
which establishes the first claim. The calculation of $M_{k,p}$ is immediate.
By the uniqueness of the Bruhat decomposition, we have
$$v_{w_1,p}=\frac{\#\left(\left(\begin{smallmatrix}1 &  &    \\ & 1  &  \F_p\\ &   & 1 \end{smallmatrix}\right)\left(\begin{smallmatrix}1 &   &    \\  &   &  1\\  & 1  & \end{smallmatrix}\right)P_{21}(\F_p)\right)}{\#\GL_3(\F_p)}=\frac{\#\F_p\#P_{21}(\F_p)}{\#\GL_3(\F_p)}=pv_1,$$
and similarly
$$v_{w_2,p}=\frac{\#\left(\left(\begin{smallmatrix}1 & \F_p & \F_p   \\  & 1  &  \\  &   & 1 \end{smallmatrix}\right)\left(\begin{smallmatrix}1 &   &    \\ &   &  1\\ & 1  & \end{smallmatrix}\right)P_{21}(\F_p)\right)}{\#\GL_3(\F_p)}=p^2v_1.$$
The result now follows since $v_1+v_2+v_3=1$.

\begin{lemma}\label{reps2}
The set $$X_p=\left\{1,
w=\left(\begin{smallmatrix}
1 &   &    \\
  &   &  1\\
  & 1  & 
\end{smallmatrix}\right)
\right\}$$ is a set of representatives of
$P_{21}(\Z_p) \bs \GL_3(\Z_p) / \Gamma_0(p) $.
The corresponding $M_{k,p}$ and volumes $v_{k,p}$ are given by
\begin{itemize} 
\item for  $k=1$ we have $M_{k,p}=M \cap K_p$, and $v_{1,p}=\frac{p-1}{p^3-1}$.
\item for $k=\left(\begin{smallmatrix}
 1&   &    \\
  &   &  1\\
  &  1 & 
\end{smallmatrix}\right)$ we have $M_{k,p}=\left(\begin{smallmatrix}
 *& *  &   \\
 p\Z_p & *  &  \\
  &   & *
\end{smallmatrix}\right)$,  and $v_{k,p}=1-v_{1,p}$.
\end{itemize}
\end{lemma}
\textbf{Proof.}
  Similar to Lemma~\ref{reps1}.

 \section{Local computations II: oldforms and orthonormal bases}
 
 In this section we prove the propositions in  Subsection  \ref{13}.\\

\textbf{Proof} of Proposition \ref{prop3}.
By Section~2 of~\cite{Ree}, the space of oldforms generated by $f$ for the group $\Gamma_0(p)$ is generated by three forms $T_0f=f$, $T_1f$, $T_2f$ 
where 
$$T_1 f(z)= \frac1{\sqrt{p}} \sum_{h \in \Delta(p)} f(\mt{h}{}{}{1} z)$$
where $\Delta(p)=\left\{ \mt{1}{k}{}{p} : 0 \le k \le p-1\right\} \cup  \mt{p}{}{}{1}$
and 
$$
T_2 f (z)= f\left(\left(\begin{smallmatrix} p &  &  \\   & p & \\  &  & 1\end{smallmatrix}\right)z\right).$$
The Fourier expansion of $f$ is given by
\begin{align*}
 f(z) = \sum_{\gamma \in \Gamma_\infty \bs \SL_2(\Z)}
    \sum_{m_1=0}^{\infty} \sum_{m_2 \neq 0}&
    \frac{A_\pi(m_1,m_2)}{|m_1m_2|} 
    e\left(m_1\pi_2(\gamma v)+m_2\Re(\gamma.z_2)\right) \\
    &\times  W_{\nu_\pi}\left(M
\begin{pmatrix}\frac{y_1y_2}{|j(\gamma,z_2)|} &  &   \\ & y_1 |j(\gamma,z_2)| &  \\    &    & 1
\end{pmatrix}\right)
\end{align*}
where $$z=xy, \quad 
x=\left(\begin{matrix}
1 & x_{12} & x_{13} \\& 1 & x_{23} \\ &  & 1\end{matrix}\right),  \quad 
y=\left(\begin{matrix}y_1y_2 &  & \\& y_1 &  \\  &  & 1
\end{matrix}\right), \quad 
z_2=x_{12}+iy_2, \quad 
v=\left(\begin{matrix}
x_{13}\\
x_{23}
\end{matrix}\right),$$
$M = \text{diag}(|m_1m_2|, m_1, 1)$, $\pi_2: \R^2 \to \R$ is the projection on the second coordinate
and $j(\gamma,z)=cz+d$ for
$\gamma=\mt{a}{b}{c}{d}$.
Note that if $h=\mt{a}{b}{}{d} \in \GL_2(\R)$, then
$\mt{h}{}{}{1}z=x'y'$,
where
$$y'=\left(\begin{matrix}y_1y_2a &  &  \\& y_1d &  \\ &   & 1
\end{matrix}\right), \quad 
x'=\left(\begin{matrix}1 &  \frac{x_{12}a +b}{d} &   x'_{13} \\  & 1 & x'_{23} \\ &  & 1 \end{matrix}\right), \quad 
 \left(\begin{matrix}
x'_{13}\\
x'_{23}
\end{matrix}\right)=hv.$$
In particular $y'_1=j(h,z_2)y_1$, $y'_2=\frac{\det(h)}{j(h,z_2)^2}y_2$ and $z'_2=hz_2$.
Using that and
 $$(\Gamma_\infty \bs \SL_2(\Z))\Delta(p)=\Delta(p)(\Gamma_\infty \bs \SL_2(\Z)),$$
 we can write $\sqrt{p}(T_1 f)(z)$ as 
 \begin{displaymath}
 \begin{split}
& \sum_{h \in \Delta(p)}\sum_{\gamma \in \Gamma_\infty \bs \SL_2(\Z)}\sum_{m_1=0}^{\infty} \sum_{m_2 \neq 0} \frac{A_\pi(m_1,m_2)}{|m_1m_2|} e\left(m_1 \pi_2(\gamma h v)+m_2\Re(\gamma(hz_2))\right) \\ &\quad\quad\quad\quad\quad \times W_{\nu_\pi}\left(
 M\left(\begin{smallmatrix}\frac{y_1y_2}{|j(\gamma,hz_2)|}\frac{p}{j(h,z_2)} &              &    \\ & y_1 |j(\gamma,hz_2)|j(h,z_2) &  \\&  & 1\end{smallmatrix}\right)\right)\\ 
 = & \sum_{h \in \Delta(p)}\sum_{\gamma \in \Gamma_\infty \bs \SL_2(\Z)} \sum_{m_1=0}^{\infty} \sum_{m_2 \neq 0} \frac{A_\pi(m_1,m_2)}{|m_1m_2|} e\left(m_1 \pi_2(\gamma h v)+m_2\Re(\gamma(hz_2))\right) \\ 
 & \quad\quad\quad\quad\quad \times W_{\nu_\pi}\left(
 M\left( \begin{smallmatrix}\frac{y_1y_2p}{|j(\gamma h,z_2)|} &   &  \\ & y_1 |j(\gamma h,z_2)| &  \\ &  & 1\end{smallmatrix}\right)\right)\\
 = & \sum_{k =0}^{p-1} \sum_{\gamma \in \Gamma_\infty \bs \SL_2(\Z)}\sum_{m_1=0}^{\infty} \sum_{m_2 \neq 0}\frac{A_\pi(m_1,m_2)}{|m_1m_2|}  e\left(m_1 p\pi_2(\gamma v)+\frac{m_2}{p}\left(\Re(\gamma z_2)+k\right)\right) \\ 
& \quad\quad\quad\quad\quad 
\times  W_{\nu_\pi}\left(
 M\left( \begin{smallmatrix} \frac{y_1y_2}{|j(\gamma ,z_2)|} &  &  \\  & p y_1 |j(\gamma,z_2)| &  \\ &  & 1\end{smallmatrix}\right)\right)\\
 +& \sum_{\gamma \in \Gamma_\infty \bs \SL_2(\Z)} \sum_{m_1=0}^{\infty} \sum_{m_2 \neq 0} \frac{A_\phi(m_1,m_2)}{|m_1m_2|}  e\left(m_1 \pi_2(\gamma v)+pm_2\Re(\gamma.z_2)\right)  W_{\nu_\pi}\left(
 M\left(\begin{smallmatrix}\frac{y_1y_2p}{|j(\gamma,z_2)|} &  &   \\& y_1 |j(\gamma,z_2)| &  \\ & & 1\end{smallmatrix}\right)\right),
 \end{split}
 \end{displaymath}
and the desired expression for the Fourier coefficients of $T_1 f$ follows. 

The calculation of the Fourier coefficients of $T_2 f$ is similar, but simpler. \medskip

The calculation of the Gram matrix is Rankin-Selberg theory.
Let $$E(z,s)=\sum_{\Gamma_0(\infty) \bs \Gamma_0(p)}
\det(\gamma z)^s,$$
where $\Gamma_0(\infty)=\left(\begin{smallmatrix}
    * & * & * \\ * & *&  *\\ & & 1\end{smallmatrix}\right) \cap \GL_3(\Z)$.
Then $E(z,s)$ has a simple pole at $s=1$ with 
    residue~$\frac{2\pi}{3\zeta(3)}\frac{p-1}{p^3-1}.$
On the other hand, by  unfolding 
we obtain for $0 \le i, j \le 2$ and $\Re(s)$ large enough
\begin{equation}\label{RankinSelberg}
    \langle T_i f E(\cdot,s), T_j f \rangle= 
\sum_{m_1,m_2>0}
\frac{A_\pi^{(i)}(m_1,m_2)\overline{A_\pi^{(j)}(m_1,m_2)}}{m_1^{2s}m_2^s}G_{\nu_\pi,\nu_\pi}(s),
\end{equation} 
where 
$$G_{\nu_1,\nu_2}(s)=\int \int W_{\nu_1}(y)\overline{W}_{\nu_2}(y) \, (y_1^2y_2)^s  d^\times y.$$
Taking the residue of~(\ref{RankinSelberg}) at $s=1$ with $i=j=0$ gives 
$$\|f\|^2=\frac3{2\pi}\frac{p^3-1}{p-1}G_{\nu_\pi,\nu_\pi}(1) \cdot\underset{s=1}{\rm{res}} L(\pi \times \tilde{\pi},s).$$
On the other hand, taking $i=0$ and $j=1$ we obtain  \begin{align*}
 \frac1{\sqrt{p}}\langle f E(\cdot,s), T_1 f \rangle
&=G_{\nu_\pi,\nu_\pi}(s)\sum_{m_1,m_2>0}
\frac{A_\pi(m_1,m_2)\big(\overline{A_{\pi}(\frac{m_1}{p},pm_2)}+\overline{A_\pi(m_1,\frac{m_2}p)}\big)}{m_1^{2s}m_2^s}\\
&=G_{\nu_\pi,\nu_\pi}(s)
\sum_{\substack{m_1,m_2>0\\ p \nmid m_1m_2}}\frac{|A_\pi(m_1,m_2)|^2}{m_1^{2s}m_2^s}\sum_{i, j \ge 0}
 \frac{A_\pi(p^i,p^j) (\overline{A_{\pi}(p^{i-1},p^{j+1})}
        +\overline{A_\pi(p^i,p^{j-1})})}{p^{(2i+j)s}}.
    \end{align*}
Using the Hecke relation
$$A_\pi(p,1)A_\pi(p^i,p^j)= A_\pi(p^{i+1},p^j)+A_\pi(p^i,p^{j-1})+A_\pi(p^{i-1},p^{j+1})$$
and the following expression for the unramified (recall that $\pi$ has level one) local $L$-factor
$$L_p(\pi \times \tilde{\pi},s )={\zeta_p(3s)}\sum_{i,j \ge 0} |A_\pi(p^i,p^j)|^2 p^{-(2i+j)s},$$
a simple calculation shows that 
$$\sum_{i, j \ge 0}
 \frac{A_\pi(p^i,p^j) (\overline{A_{\pi}(p^{i-1},p^{j+1})}
        +\overline{A_\pi(p^i,p^{j-1})})}{p^{(2i+j)s}}=(p^{-s}-p^{-3s}){A_\pi(1,p)}L_p(\pi \times \tilde{\pi},s),$$
and thus taking the residue at $s=1$ we obtain
the desired expression for $ \langle f , T_1 f \rangle$.

The inner product  $\langle f , T_2 f \rangle$ is calculated in a similar way. 

Changing variables in the definition of the inner product, we have $\langle T_1 f, T_2 f \rangle = \langle f, T_1 f \rangle$ and $\|T_2f\|^2=\|f\|^2$.

To compute $\|T_1f\|^2$, let $g=\sqrt{p}A_\pi(p,1)f-T_1f.$
Using the Hecke relations, the Fourier coefficients of $g$ are given by $\sqrt{p}A_\pi(pm_1,m_2)$.
    On the other hand
    $$\langle T_1 f+g, T_1 f-g \rangle = \|T_1f\|^2 - \|g\|^2 +2i\Im(\langle g, T_1 f \rangle),$$
    and thus $$\|T_1 f\|^2=\|g\|^2 + \sqrt{p} \Re (A_\pi(p,1)\langle f, T_1f-g\rangle). $$
     By Rankin-Selberg unfolding we find
    $$\langle f,g  \rangle=p^{3/2}  \langle T_2 f,f  \rangle,$$
    from which follows 
    $$\sqrt{p} A_\pi(p,1)\langle f, T_1f-g\rangle=
    -\frac{p^4-2p^3+p}{p^3-1}|A(1,p)|^2\|f\|^2.$$
   Finally, using Rankin-Selberg unfolding once more, we have
    $$\| g \|^2 = p^3\underset{s=1}{{\rm{res}}}\left(1
    -\zeta_p(3s)\frac{\mathcal{L}_p(\pi,s)}{L_p(\pi \times \tilde{\pi},s)}  \right)\|f\|^2,
    $$
    where we have defined $$\mathcal{L}_p(\pi,s)=\sum_{j \ge 0} \frac{|A_\pi(1,p^j)|^2}{p^{js}}.$$
    To compute $\mathcal{L}_p(\pi,s)$, 
   we express the Fourier coefficient $A_\pi(1,p^j)$ in terms of the Satake parameters using the Schur polynomials and we 
compute the sum of the geometric series. Using $\alpha_p\beta_p\gamma_p=1$, we find
\begin{align*}
\frac{\mathcal{L}_p(\pi,s)}{L_p(\pi \times \tilde{\pi},s)}
= 1& -(\alpha_p\beta_p+\beta_p\gamma_p+\gamma_p\alpha_p)
(\overline{\alpha_p\beta_p}+\overline{\beta_p\gamma_p}+\overline{\gamma_p\alpha_p})p^{-2s}\\
& +
((\alpha_p+\beta_p)(\beta_p+\gamma_p)(\gamma_p+\alpha_p)+
(\overline{\alpha_p}+\overline{\beta_p})(\overline{\beta_p}+\overline{\gamma_p})(\overline{\gamma_p}+\overline{\alpha_p}))p^{-3s}\\
& - (\alpha_p+\beta_p+\gamma_p)(\overline{\alpha_p}+\overline{\beta_p}+\overline{\gamma_p})p^{-4s} +p^{-6s}\\
= & \quad(1-p^{-3s})^2-|A_\pi(1,p)|^2p^{-2s}(1-p^{-s})^2.
\end{align*}
Taking the residue at $s=1$ and simplifying finishes the proof.\\

We now turn to the continuous spectrum. 
Since the maximal Eisenstein series are the most complicated case, 
we present it first.\\

\textbf{Proof} of Proposition \ref{prop5}. 
Let $P =P_{21}$ and write its Levi decomposition $P=NM$.
 Since $M \cong \GL_2 \times \GL_1$ we know $\sigma$ is of the form
    $\tau \otimes \chi$ where $\tau$ is an automorphic representation of $\GL_2$ and $\chi$ is a character.
    Since $\I_P(\sigma_s)^{\Gamma^{(3)}_0(p)}$ is non-trivial, Lemma~\ref{Isom} implies that $\sigma^{M_k}$ is non-trivial for some $k$.
    The computation of $M_k$ in Lemma~\ref{reps2} implies that 
    $\tau$ must have a non-zero $\SO_2(\R)\Gamma^{(2)}_0(p)$-fixed vector and $\chi$ is unramified at $p$. 
    Moreover, $\chi$ is also unramified at every other finite places and $\chi$
    is trivial on $\R^+$.
    So either $\chi=1$ or $\chi=\sgn$.
    Finally, the central character of $\I_P(\sigma_s)$ is 
    $\omega_\tau \chi$, hence $\omega_\tau=\chi^{-1}$.
    But since $\mt{-1}{}{}{-1} \in \SO_2(\R)$, 
    we must have $\omega_\sigma=\chi=1.$
    We have two cases to consider:

    a) If $a(\tau)=1$ we have nothing to prove.

    b) If $a(\tau)=0$, let $\phi$ be the spherical vector in
    $\tau$. By abuse of notation, also denote by $\phi$ the spherical vector in $\I_P(\sigma_s)$, it is given by
    $$\phi(nmk)=\phi(\mathfrak{m}_2(m)).$$
     The space of oldforms generated by
$E(\cdot,\phi,s)$ is generated by 
$T_0 E(\cdot,\phi,s),$ $T_1 E(\cdot,\phi,s)$ and 
$T_2 E(\cdot,\phi,s)$, where $T_0,T_1,T_2$ are as in the cuspidal case. 
Let $X=\{1,w\}$ where $w$ is the element of
$\GL_3(\A)$ whose component at each place $v\neq p$ equals $1$ and whose component at $p$ is given by the corresponding $w$ in Lemma~\ref{reps2}. Then $X$ is a set of representatives of 
$(P \cap K) \bs K / \Gamma_0(p)$ and correspondingly we  
consider the isomorphism given by~(\ref{isomorphism}).
For convenience, we use the notation 
$\phi_{|p}=\tau\left(\mt{p^{-1}}{}{}{1}\right)\phi.$

\begin{lemma}\label{Ttilde1}
We have
   $T_1E(\cdot,\phi,s)= E(\cdot,\tilde{T}_1(s)\phi,s)$ and $T_2E(\cdot,\phi,s)= E(\cdot,\tilde{T}_2(s)\phi,s)$ where the images of $\tilde{T}_1(s)\phi$
and $\tilde{T}_2(s) \phi$ under isomorphism~(\ref{isomorphism}) are respectively given by
\begin{equation*}
    \begin{split}
    (\tilde{T}_1(s) \phi)_1(m)&=p^{\frac12}(\alpha_p+\beta_p)
    \phi(\mathfrak{m}_2(m)),\\
    (\tilde{T}_1(s) \phi)_w(m)&=
    p^{-1/2-2s}
    \phi(\mathfrak{m}_2(m))
    +
    p^{s}
    \phi_{|p}\left(\mathfrak{m}_2(m)\right),
    \end{split}
    \quad
    \begin{split}
        (\tilde{T}_2(s) \phi)_1(m)&=p^{1+2s}
    \phi(\mathfrak{m}_2(m)),\\
    (\tilde{T}_2(s) \phi)_w(m)&=
    p^{-\frac12-s} 
    \phi_{|p}(\mathfrak{m}_2(m)).
    \end{split}
\end{equation*}
\end{lemma}

\textbf{Proof.}
    We focus on $T_1E(\cdot,\phi,s)$ since $T_2$ is simpler.
       Using~(\ref{Intertwin}) it suffices to show that
    $$\frac1{\sqrt{p}} \sum_{h \in   \Delta(p)} \I_P\left(\tau_s,\mt{h^{-1}}{}{}{1}\right)\phi=\tilde{T_1}(s)\phi.$$
So we need to show that the left hand side and the right hand side coincide on 
$N(\A)M(\A) \Gamma_0(p)$ and on $N(\A)M(\A) w \Gamma_0(p)$.
Let $\gamma=\left(\begin{smallmatrix}
A &  \begin{smallmatrix} \gamma_1  \\ \gamma_2 \end{smallmatrix} \\ \begin{smallmatrix} p \gamma_4 & p\gamma_5 \end{smallmatrix} & \gamma_3
\end{smallmatrix}\right) \in \Gamma_0(p)$
and let $h \in \Delta(p)$.
We have
\begin{equation*}\label{decomposition}
    \gamma \left(\begin{matrix}
h^{-1} &  \begin{matrix} \phantom{\gamma_1}  \\ \phantom{\gamma_2} \end{matrix} \\ \begin{matrix} \phantom{\gamma_4'} & \phantom{\gamma_5'} \end{matrix} & 1\end{matrix}\right) = \left(\begin{matrix}
Ah^{-1} &  \begin{matrix} \gamma_1  \\ \gamma_2 \end{matrix} \\ \begin{matrix} \gamma_4' & \gamma_5' \end{matrix} & \gamma_3\end{matrix}\right)
=\left(\begin{matrix}Ah^{-1} &  \begin{matrix}  \phantom{\gamma_1} \\ \phantom{\gamma_2} \end{matrix} \\
 \begin{matrix}  &  \end{matrix} & 1 \end{matrix}\right)\left(\begin{matrix}1 & & \gamma_1'\\ & 1 & \gamma_2'\\ \gamma'_4 & \gamma'_5 & \gamma_3
\end{matrix}\right)
\end{equation*}
where
$( \begin{matrix} \gamma_4' & \gamma_5' \end{matrix})=
 (\begin{matrix} p\gamma_4 & p\gamma_5 \end{matrix}) h^{-1} \in \Z_p^2$ 
 and $\left(\begin{matrix}  \gamma'_1 \\ \gamma'_2 \end{matrix} \right)=hA^{-1}\left(\begin{matrix}  \gamma'_1 \\ \gamma'_2 \end{matrix} \right) \in \Z_p^2$.
 In particular, $\left(\begin{smallmatrix}
1 & & \gamma_1'\\
 & 1 & \gamma_2'\\
\gamma'_4 & \gamma'_5 & \gamma_3
\end{smallmatrix}\right) \in K_p$
and thus 
$I_{\mathbf{s}_P}(\gamma h)=p^{\frac12+s}.$
It follows that for all $n \in N(\A)$, $m\in M(\A)$ we have
\begin{align*}
    \I_P(\tau_s,\mt{h^{-1}}{}{}{1})\phi(nm\gamma)&=p^{\frac12+s}\phi\left(nm\left(\begin{smallmatrix}
Ah^{-1} &  \begin{smallmatrix}  \phantom{\gamma_1} \\ \phantom{\gamma_2} \end{smallmatrix} \\
 \begin{smallmatrix}  &  \end{smallmatrix} & 1
\end{smallmatrix}\right)\right),
\end{align*}
and thus by definition of the $\GL_2$ Hecke operators
 $$\frac1{\sqrt{p}} \sum_{h \in   \Delta(p)} \I_P\left(\tau_s,\mt{h^{-1}}{}{}{1}\right)\phi(nm\gamma)=p^{\frac12+s} 
\lambda_\phi(p) 
\phi(\mathfrak{m}_2(m)),$$
where $\lambda_\phi(p)=a_p+b_p$ is the Hecke eigenvalue of $\phi$.
Next write $A=\mt{a}{b}{c}{d} \in \GL_2(\Z_p)$. We have
\begin{align*}
\left(\begin{matrix}a & b \\c &d \end{matrix}\right)
\left(\begin{matrix}p^{-1} &  \\ &1 \end{matrix}\right)
&=
\left(\begin{matrix} p^{-1}a &b\\ p^{-1}c& d\end{matrix}\right)=\begin{cases}
\left(\begin{matrix}p^{-1}& \\ & 1 \end{matrix}\right)
\left(\begin{matrix} a& pb \\ p^{-1}c & d\end{matrix}\right), &p \mid c,\\
\left(\begin{matrix} 1 &p^{-1}ac^{-1}\\  & p^{-1}\end{matrix}\right)
\left(\begin{matrix} & b-adc^{-1} \\ c & pd\\ \end{matrix}\right), & p \nmid c,
\end{cases}
\end{align*}
and 
\begin{align*}
\left(\begin{matrix} a &b \\c &d\end{matrix}\right)&
\left(\begin{matrix} 1 &-p^{-1}k\\ & p^{-1}\end{matrix}\right)
=\left(\begin{matrix} a & p^{-1}(b-ka) \\c & p^{-1}(d-kc)\end{matrix}\right)\\
&=\begin{cases}
\left(\begin{matrix} p^{-1} & \\ &1\end{matrix}\right)
\left(\begin{matrix}pa &b-ka\\c &p^{-1}(d-kc)\end{matrix}\right), &p \mid d-kc,\\
\left(\begin{matrix}1 &p^{-1}(b-ka)(d-kc)^{-1}\\& p^{-1}\end{matrix}\right)
\left(\begin{matrix} a-c(b-ka)(d-kc)^{-1}& \\pc &d-kc\end{matrix}\right), & p \nmid d-kc.
\end{cases}
\end{align*}
Now note that either $p \mid c$ in which case $p \nmid d-kc$ for all
integer $k$, or $p \nmid c$ and there exists a unique $k \pmod c$ such that $p \mid d-kc$.
This implies that for any fixed $\gamma \in \Gamma_0(p)$, there exists
a unique $h_0 \in \Delta(p)$ such that
$w\gamma \mt{h_0^{-1}}{}{}{1} \in \left( \begin{smallmatrix}
    p^{-1}& & \\ & 1 & \\ & & 1\end{smallmatrix}\right)K_p$, 
and for the remaining $p$ elements $h \neq h_0$ in $\Delta (p)$ we
have $w\gamma \mt{h^{-1}}{}{}{1} \in N(\Q_p)\left(\begin{smallmatrix}
    1& & \\ & 1 & \\ & & p^{-1}\end{smallmatrix}\right)K_p$.
Then for all $n \in N(\A)$ and $m \in M(\A)$ we have 
$$\I_P\left(\tau_s,\mt{h^{-1}}{}{}{1}\right)\phi(nmw\gamma)=\begin{cases}
    p^{\frac12+s}\phi\left(m\left(\begin{smallmatrix}p^{-1} & & \\ & 1& \\ & & 1\end{smallmatrix}\right)\right), & h=h_0\\
    p^{-1-2s}\phi\left(m\left(\begin{smallmatrix}1 & & \\ & 1& \\ & & p^{-1}\end{smallmatrix}\right)\right), & h \neq h_0.
\end{cases}$$
Hence
 $$\frac1{\sqrt{p}} \sum_{h \in   \Delta(p)} \I_P\left(\tau_s,\mt{h^{-1}}{}{}{1}\right)\phi(nmw\gamma)=p^{s} 
\phi_{|p}(\mathfrak{m}_2(m))+p^{-\frac12-2s}\phi(\mathfrak{m}_2(m)),$$
which establishes the result.
The calculation for $T_2$ is similar.\\

Using $\GL_2$ Rankin-Selberg theory and Lemma~\ref{reps2} one can now compute the Gram matrix of $\phi, \tilde{T}_1(s)\phi$ and $\tilde{T}_2(s)\phi$ with respect to the inner product~(\ref{InnerProduct}).
As mentioned in the introduction, it turns out to be identical as in the cuspidal case.\\

\textbf{Proof} of Proposition \ref{prop4}. For $P=P_{111}$ we have $M=A$ and thus the representation~$\sigma$~is a character
of~$A(\Q)A(\R)^+\bs A(\A)$ that is trivial on the centre.
By Lemmas~\ref{Isom} and~{\ref{reps1}},~$\sigma$~is also trivial on $A \cap K$ and thus $\sigma=1$. For any given $s=(s_1,s_2) \in (i\R)^2$, the space of oldforms generated by
$E(\cdot,1,s)$ is generated by 
$T_0 E(\cdot,1,s),$ $T_1 E(\cdot,1,s)$ and 
$T_2 E(\cdot,1,s)$, where $T_0,T_1,T_2$ are as in the cuspidal case. 
Let $X=\{1,w_1,w_2\}$ where $w_i$ is the element of
$\GL_3(\A)$ whose component at each place $v\neq p$ equals $1$ and whose component at $p$ is given by the corresponding $w_i$ in Lemma~\ref{reps1}. The following is the analogue of Lemma~\ref{Ttilde1}. 
\begin{lemma}
    We have  $T_1 E(\cdot,1,s)=E(\cdot,\tilde{T_1}(s)1,s)$
    and $T_2 E(x,1,s)=E(x,\tilde{T}_2(s) 1, s)$
    where the images under~(\ref{isomorphism}) of $\tilde{T_1}(s)1$
    and $\tilde{T}_2(s) 1$ are respectively given by
    \begin{equation*}
        \begin{split}
            (\tilde{T_1}(s)1)_1&=p^{1/2}(\alpha_p+\beta_p),\\
            (\tilde{T_1}(s)1)_{w_1}&=p^{1/2}\alpha_p+p^{-1/2}\gamma_p,  \\
            (\tilde{T_1}(s)1)_{w_2}&=p^{-1/2}(\beta_p+\gamma_p),
        \end{split}
        \quad
        \begin{split}
            (\tilde{T}_2(s) 1)_1&=p \gamma_p^{-1},\\
            (\tilde{T}_2(s) 1)_{w_1}&=\beta_p^{-1},\\
            (\tilde{T}_2(s) 1)_{w_2}&=p^{-1}\alpha_p^{-1}.
        \end{split}
    \end{equation*}
\end{lemma}
\textbf{Proof.}
    Similar to proof of Lemma~\ref{Ttilde1}.\\

Using Lemma~\ref{reps1} it is now easy to calculate the Gram
matrix of  $1, \tilde{T}_1(s)1$ and $\tilde{T}_2(s)1$ with respect to the inner product~(\ref{InnerProduct}). Again, the result is identical to the
cuspidal case.\\

\section{The Kuznetsov formula}\label{seckuz}

We quote the ${\rm GL}(3)$ Kuznetsov formula as presented in \cite[Sections 2.3 - 2.6]{BL} but with Kloosterman sums for general level imported from  \cite[Theorem 6]{BBM}. 

\subsection{Test functions and integral transforms} 
Our test functions $h$ live on the two-dimensional hyperplane in $\Bbb{C}^3$ of triples $\mu = (\mu_1, \mu_2, \mu_3)$ satisfying $\mu_1+\mu_2 + \mu_3 = 0$. We recall that $\mu$ occuring as Langlands parameters of an automorphic form satisfies 
\begin{equation}\label{auto}
   |\Re \mu_j| \leq 5/14, \quad \{\mu_1, \mu_2, \mu_3\} = \{-\overline{\mu}_1, -\overline{\mu}_2, -\overline{\mu}_3\}
   \end{equation}
   by the Kim-Sarnak bound and unitarity.

\begin{defi}\label{def1} For $ A_0 > 0$ let  $\tilde{\mathcal{H}}(A_0)$   be the class  of test functions $h$ that are Weyl group invariant, holomorphic for $|\Re \mu_j | \leq A_0$, bounded by $O_A(\| \mu \|^{-A})$ in this strip for every $A > 0$ 
 and have zeros at 
\begin{equation}\label{zeros}
\begin{split}
   & \mu_i - \mu_j = n \equiv 1 \, (\text{\rm mod } 2), \quad |n| \leq A_0, \quad 1 \leq i < j \leq 3. 
    \end{split}
    \end{equation}
    We define $\mathcal{H}(A_0) \subseteq \tilde{\mathcal{H}}(A_0)$ to be the subset of those functions that have in addition zeros at
    \begin{equation}\label{zeros1}
\begin{split}
  &\tfrac{1}{2} +\mu_i + \mu_j = n \equiv 0\, (\text{\rm mod } 2), \quad |n| \leq A_0, \quad  1 \leq i \leq j \leq 3.
    \end{split}
    \end{equation}
    \end{defi}
\textbf{Remarks.} A typical function we have in mind is
\begin{equation}\label{testh}
\begin{split}
h(\mu) =-\exp(\mu_1^2 +\mu_2^2 + \mu_3^2)& \Big(\prod_{\substack{|n| \leq A_0\\ n \equiv 1 \, (\text{mod }2)}}\prod_{1 \leq i < j \leq 3}(\mu_i - \mu_j - n)\Big) \\
&\times \Big(\prod_{\pm} \prod_{\substack{|n| \leq A_0\\ n \equiv 0 \, (\text{mod }2)}} \prod_{1 \leq i \leq j \leq 3}(\tfrac{1}{2} \pm \mu_i \pm \mu_j - n )\Big).
\end{split}
\end{equation}
Note that this function is strictly positive on spectral parameters satisfying \eqref{auto}. 
 
 Also note that the holomorphicity of $h$ and the upper bound  $h(\mu) \ll \| \mu \|^{-A}$ imply bounds for all derivatives.  The   zeros in \eqref{zeros} are needed to optimize the performance of the Kuznetsov formula, the zeros in \eqref{zeros1} are needed to obtain a suitable approximate functional equation.\\


We define the integral transforms
\begin{align*}
    \begin{aligned}
& \Phi_{w_4}(y)=\int_{\Re \mu=0} h(\mu) K_{w_4}(y ; \mu) \operatorname{spec}(\mu) \mathrm{d} \mu, \\
& \Phi_{w_5}(y)=\int_{\Re \mu=0} h(\mu) K_{w_4}(-y ;-\mu) \operatorname{spec}(\mu) \mathrm{d} \mu, \\
& \Phi_{w_6}(y_1, y_2)=\int_{\Re \mu=0} h(\mu) K^{\text{sym}}_{w_6}\left(\left(y_1, y_2\right) ; \mu\right) \operatorname{spec}(\mu) \mathrm{d} \mu .
\end{aligned}
\end{align*}
where
$$\text{spec}(\mu) = (\mu_1 - \mu_2)(\mu_2 - \mu_3)(\mu_3 - \mu_1)  \tan \Big(\frac{\pi}{2}(\mu_1 - \mu_2)\Big)\tan \Big(\frac{\pi}{2}(\mu_2 - \mu_3)\Big)\tan \Big(\frac{\pi}{2}(\mu_3 - \mu_1)\Big)$$
and
\begin{displaymath}
\begin{split}
    K_{w_4}(y;\mu)&=\int_{-i\infty}^{i\infty}|y|^{-s}\tilde{G}^{\mathrm{sgn}(y)}(s,\mu)\frac{\mathrm{d}s}{2\pi i},\\
     K_{w_6}^{\mathrm{sym}}(y ; \mu)&=\int_{-i\infty}^{i\infty} \int_{-i\infty}^{i\infty} \left(\pi^2 |y_1|\right)^{-s_1}\left(\pi^2 |y_2|\right)^{-s_2} G_{\mathrm{sym}}^{\mathrm{sgn}(y_1),\mathrm{sgn}(y_2)}(s, \mu) \frac{\mathrm{d} s_1 \mathrm{~d} s_2}{(2 \pi i)^2}
     \end{split}
\end{displaymath}    
with
\begin{displaymath}
\begin{split}
    \tilde{G}^{ \pm}(s, \mu)&:=\frac{\pi^{-3 s}}{12288 \pi^{7 / 2}}\left(\prod_{j=1}^3 \frac{\Gamma\left(\frac{1}{2}\left(s-\mu_j\right)\right)}{\Gamma\left(\frac{1}{2}\left(1-s+\mu_j\right)\right)} \pm i \prod_{j=1}^3 \frac{\Gamma\left(\frac{1}{2}\left(1+s-\mu_j\right)\right)}{\Gamma\left(\frac{1}{2}\left(2-s+\mu_j\right)\right)}\right)\\
   G_{\mathrm{sym}}^{\epsilon_1,\epsilon_2}(s,\mu)&=\frac{1}{1024 \pi^{5 / 2}} \sum_{\substack{d_1,d_2,d_3 \in\{0,1\}\\ d_3\equiv d_1+d_2 \, (\text{mod } 2)}} \epsilon_1^{d_1} \epsilon_2^{d_2}(-1)^{d_1 d_2} \frac{\Gamma (\frac{1+d_3-s_1-s_2}{2} )}{\Gamma (\frac{d_3+s_1+s_2}{2} )} \prod_{i=1}^3 \frac{\Gamma (\frac{d_1+s_1-\mu_i}{2} ) \Gamma (\frac{d_2+s_2+\mu_i}{2} )}{\Gamma(\frac{1+d_1-s_1+\mu_i}{2} ) \Gamma (\frac{1+d_2-s_2-\mu_i}{2} )}.
   \end{split}
\end{displaymath}
The precise shape of the integral kernels is not relevant for the subsequent discussion and only included for completeness. All we need is summarized in the following lemma which is a restatement of \cite[Lemma 2.4, 2.5, 2.6]{BL}.  

\begin{lemma}\label{lem.Phi4Truncation} 
     Let $A > 0$,   $j, j_1, j_2 \in \Bbb{N}_0$. There is a choice of  $A_0 $  such that for $h \in \mathcal{H}(A_0)$ we have  
     \begin{equation}\label{lem8a}
       |y|^j \Phi^{(j)}_{w_4}(y), \,\, |y|^j\Phi^{(j)}_{w_5}(y)\ll    |y|^{j/3} \min\big(|y|^A, |y|^{-1/6+\varepsilon}\big),\\
    \end{equation}
      \begin{equation}\label{lem8b}
      \begin{split}
       y_1^{j_1}y_2^{j_2}&\frac{\mathrm{d}^{j_1+j_2}}{\mathrm{d} y_1^{j_1}\mathrm{d} y_2^{j_2}}\Phi_{w_6}(y_1,y_2)\\
       &\ll  
         (1+|y_1|^{1/2}+|y_1|^{1/3}|y_2|^{1/6})^{j_1}(1+|y_2|^{1/2}+|y_2|^{1/3}|y_1|^{1/6})^{j_2}\min\big(|y_1|,|y_2|\big)^A,
         \end{split}
\end{equation}
          \begin{equation}\label{lem8c}
      \Phi_{w_6}(y_1,y_2)\ll  
      (2 + |\log |y_1|| + |\log |y_2||)^3(1 + |y_1|)^{-1/4}(1 + |y_2|)^{-1/4}, 
      \end{equation}
     with   implicit constants depending on $h, A,  j, j_1, j_2$. 
   \end{lemma}

\subsection{Kloosterman sums} 
For   integers $n_1,n_2,m_1,m_2$ and positive integers $N, D_1,D_2$, we define two types of Kloosterman sums as follows: \begin{align*}
\tilde{S}\left(n_1, n_2, m_1 ; D_1, D_2\right):=\sum_{\substack{C_1\, ( \text{mod }  D_1), C_2 \, (\text{mod } D_2) \\\left(C_1, D_1\right)=\left(C_2, D_2 / D_1\right)=1}} e\left(n_2 \frac{\bar{C}_1 C_2}{D_1}+m_1 \frac{\bar{C}_2}{D_2 / D_1}+n_1 \frac{C_1}{D_1}\right)
\end{align*}
for $D_1 \mid D_2$, and \begin{displaymath}
\begin{split}
& S^{(N)}\left(n_1, n_2, m_1, m_2 ; D_1, D_2\right) 
\\& =\sum_{\substack{B_1, C_1 \, (\text{mod }  D_1) \\
B_2, C_2\, (\text{mod } D_2) \\
D_1 C_2+B_1 B_2+D_2 C_1 \equiv 0 \, (\text{mod }D_1 D_2) \\
\left(B_j, C_j, D_j\right)=1, N \mid B_1}} e\left(\frac{n_1 B_1+m_1\left(Y_1 D_2-Z_1 B_2\right)}{D_1}+\frac{n_2 B_2+m_2\left(Y_2 D_1-Z_2 B_1\right)}{D_2}\right), 
\end{split}
\end{displaymath}
where $Y_j B_j+Z_j C_j \equiv 1\left(\bmod D_j\right)$ for $j=1,2$. We quote \cite[Lemma 6]{BBM}, \cite[Thm.\ 5.10]{KN} and \cite[Appendix]{BFG}. 

\begin{lemma}\label{lem1} 
{\rm (a)} If   $(t_1t_2, u_1u_2) = 1$ for $j=1, 2$ and $N \mid t_1u_1, t_2u_2$, then
\begin{displaymath}
\begin{split}
&S^{ (N)}(n_1, n_2, m_1, m_2; t_1u_1, t_2u_2) \\
&= S^{(\text{{\rm gcd}}(N, t_1))}(\bar{u}_1^2 u_2 n_1, \bar{u}_2^2 u_1 n_2, m_1, m_2; t_1,t_2) S^{({\rm gcd}(N, u_1))}(\bar{t}_1^2t_2n_1, \bar{t}_2^2t_1n_2, m_1, m_2; u_1, u_2).
\end{split}
\end{displaymath}
{\rm (b)} Let $p$ be prime and let $r_q(n)$ denote the Ramanujan sum. Then
$$S^{(p)}(n_1, n_2, m_1, m_2; p, p)= p-1 + r_p(m_1) r_p(n_2) = \begin{cases}p(p-1), & p \mid m_1, p \mid n_2,\\ p , & p  \nmid m_1n_2,\\ 0 \, & \text{else.} \end{cases}$$ 
{\rm (c)} We have 
\begin{align*}
    &S^{(1)}(n_1,n_2,m_1,m_2;D_1,D_2)\nonumber\\
    &=\sum_{D_0\mid (D_1,D_2)}D_0 \hspace{-0.3cm}\sum_{\substack{\alpha \, (\text{{\rm mod }}D_0)\\m_1\frac{D_2}{D_0}+n_2\frac{D_1}{D_0}\alpha\equiv 0 \, (\text{{\rm mod }}D_0)}} \hspace{-0.3cm}S\Big(n_1,\frac{m_1D_2+n_2D_1\alpha}{D_0^2};\frac{D_1}{D_0}\Big)S\Big(m_2,\frac{m_1D_2\overline{\alpha}+n_2D_1}{D_0^2};\frac{D_2}{D_0}\Big).
\end{align*}
{\rm (d)} We have
$$\tilde{S}(n_1, n_2, m_1; D_1, D_2) \ll \big((m_1, D_2/D_1)D_1^2, (n_1, n_2, D_1)D_2\big)(D_1D_2)^{\varepsilon}.$$
\end{lemma}

We complement this with two explicit evaluations in special cases (which is essentially \cite[Property 4.10]{BFG}). 
\begin{lemma}\label{klooster-lemma}
{\rm (a)} We have 
$S^{(p)}(n_1, n_2, m_1, m_2; p, p^2) = r_{p}(m_1)S(n_2, m_2p, p^2) .$\\
{\rm (b)}  We have
$S^{(p)}(n_1, n_2, m_1, m_2; p^2, p) =  r_{p}(m_2)S(n_1, m_1p, p^2) $. 
\end{lemma}

\textbf{Proof.} It suffices to prove (a), part (b) then follows by symmetry. To prove (a), we apply the definition and see immediately that $B_1 = 0$ in the summation condition. Hence the congruence becomes $p^3 \mid p C_2 + p^2 C_1$ which implies $p \mid C_2$, say $C_2 = p C_2'$ hence $(B_2, p) = 1$. Now the congruence reads $p \mid C_2' + C_1$ which determines $C_2'$ modulo $p$. 
We conclude
$$S^{(p)}(n_1, n_2, m_1, m_2; p, p^2) = \sum_{\substack{C_1\, (\text{mod } p)\\ (C_1, p) = 1}}  \sum_{\substack{B_2\, (\text{mod } p^2)\\ (B_2, p) = 1}} e\Big(\frac{-m_1Z_1B_2}{p} + \frac{n_2B_2}{p^2} + \frac{m_2Y_2}{p}\Big).$$
In the present situation we have $Z_1 = \bar{C}_1$ (mod $p$)   and $Y_2 \equiv \bar{B}_2$ (mod $p$), and the desired formula follows after evaluating  the $C_1$-sum as a Ramanujan sum and the $B_2$-sum as a Kloosterman sum. 

  \subsection{Spectral decomposition}

We can now state the Kuznetsov formula for the group $\Gamma_0(p) \subseteq {\rm SL}_3(\Bbb{Z})$ for $p$ prime. 

 \begin{lemma}\label{lem.Kuznetsov}
    Let $m_1,m_2,n_1,n_2\in\N$ and let $h \in \tilde{\mathcal{H}}(10)$. Let $p$ be a prime. There are absolute constants $c_{\text{{\rm cusp}}}, c_{\text{{\rm min}}}, c_{\text{{\rm max}}} > 0$ such that  \begin{equation}\label{spec}
        \mathcal{C}_{\rm new}+\mathcal{C}_{\rm old}+E_{\min}+E_{\max, {\rm new}} + E_{\max, {\rm old}} =\Delta+S_4+S_5+S_6,
    \end{equation}
    where $$
\begin{aligned}
& \mathcal{C}_{\rm new}=c_{\text{{\rm cusp}}}\sum_{\pi \in \mathcal{B}(p)} \frac{h\left(\mu_\pi\right)}{L(1,  \pi, \text{{\rm Ad}})} \overline{A_{\pi}\left(m_1, m_2\right)}A_{\pi}\left(n_1, n_2\right), \\
& \mathcal{C}_{\rm old}=c_{\text{{\rm cusp}}}\sum_{\pi \in \mathcal{B}(1)}\sum_{j=0}^2 \frac{h\left(\mu_\pi\right)}{L(1,  \pi, \text{{\rm Ad}})} \overline{B^{(j)}_{\pi}\left(m_1, m_2\right)}B^{(j)}_{\pi}\left(n_1, n_2\right), \\
& E_{\min }=c_{\text{{\rm min}}} \sum_{j=0}^2\mathop{\int\int}_{\Re s =0} \frac{h(s)}{\displaystyle\prod_{1 \leq i < l \leq 3}|\zeta(1 + s_i - s_l)|^2} \overline{B^{(j)}_s\left(m_1, m_2\right)}B^{(j)}_s\left(n_1, n_2\right) \mathrm{d} s_1 \mathrm{d}s_2, \quad (s_3 = - s_1 - s_2) \\
& E_{\max, {\rm new} }=c_{\text{{\rm max}}} \sum_{\phi \text{ {\rm level} } p} \int_{\Re s=0} \frac{h (s+t_{\phi}, s-t_{\phi},-2 s )}{L(1,  \phi, \text{{\rm Ad}})|L(1 + 3s, \phi)|^2} \overline{A_{(\phi, s)}\left(m_1, m_2\right)} A_{(\phi, s)}\left(n_1, n_2\right) \mathrm{d} s,\\
& E_{\max, {\rm old} }=c_{\text{{\rm max}}} \sum_{\phi \text{ {\rm level} } 1}\sum_{j=0}^2 \int_{\Re s=0} \frac{h (s+t_{\phi}, s-t_{\phi},-2 s )}{L(1,  \phi, \text{{\rm Ad}})|L(1 + 3s, \phi)|^2} \overline{B^{(j)}_{(\phi, s)}\left(m_1, m_2\right)} B^{(j)}_{(\phi, s)}\left(n_1, n_2\right) \mathrm{d} s,\\
\end{aligned}$$
and
$$\begin{aligned}
& \Delta=P\delta_{m_1=n_1, m_2=n_2} \int_{\Re \mu =0} h(\mu) \operatorname{spec}(\mu) \mathrm{d} \mu,\\
& S_4=P\sum_{\epsilon= \pm 1} \sum_{\substack{pD_2 \mid D_1 \\
m_2 D_1=n_1 D_2^2}} \frac{\tilde{S}\left(-\epsilon n_2, m_2, m_1; D_2, D_1\right)}{D_1 D_2} \Phi_{w_4}\left(\frac{\epsilon m_1 m_2 n_2}{D_1 D_2}\right), \\
& S_5=P\sum_{\epsilon= \pm 1} \sum_{\substack{p\mid D_1 \mid D_2 \\
m_1 D_2=n_2 D_1^2}} \frac{\tilde{S}\left(\epsilon n_1, m_1, m_2; D_1, D_2\right)}{D_1 D_2} \Phi_{w_5}\left(\frac{\epsilon n_1 m_1 m_2}{D_1 D_2}\right), \\
& S_6=P\sum_{\epsilon_1, \epsilon_2= \pm 1} \sum_{p\mid D_1, p \mid D_2} \frac{S^{(p)}\left(\epsilon_2 n_2, \epsilon_1 n_1, m_1, m_2 ; D_1, D_2\right)}{D_1 D_2} \Phi_{w_6}\left(-\frac{\epsilon_2 m_1 n_2 D_2}{D_1^2},-\frac{\epsilon_1 m_2 n_1 D_1}{D_2^2}\right). 
\end{aligned}$$
where $P = p^2 + p + 1 = [{\rm SL}_3(\Bbb{Z}) : \Gamma_0(p)]$. 
\end{lemma}
 Note that unlike in \cite{BBM} we use the probability measure on $L^2(\Gamma_0(p)\bs {\rm SL}_3(\Bbb{R})/{\rm SO}(3))$ which scales the Fourier coefficients up by a factor $P^{1/2}$ and explains the extra factor $P$ on the right hand side. 
 
\section{Proof of Theorem \ref{thm1}}

To simplify the presentation, we make the following convention. We call a term ``negligible'' if it is bounded by $\ll_A p^{-A}$ for any $A > 0$, with polynomial dependence on spectral parameters where relevant. We introduce the notation $$X \preccurlyeq Y \quad  \Longleftrightarrow \quad X \ll_{\varepsilon} Yp^{\varepsilon}.$$

\subsection{An approximate functional equation}\label{51} We start with an approximate functional equation (cf.\ \cite[Section 5]{IK}), using the information of Propositions \ref{prop1} and \ref{prop2}. Our approximate functional equation  features  two terms of \emph{different} lengths and is carefully designed to take care of the fact that the local factors at 2 and $\infty$ of the symmetric square $L$-functions can have unwanted poles. 

With this in mind we choose a small parameter
$$0 < \alpha < 1/2$$ 
to be specified later,   define
$$G(u, \mu) = (1 - 4u^2)\prod_{\pm} \prod_{\substack{|n| \leq A_0\\ n \equiv 0 \, (\text{mod } 2)}} \prod_{1 \leq i \leq j \leq 3} \big(\tfrac{1}{2}-  u + \mu_{\pi, i} + \mu_{\pi, j} + n) $$
 and consider the integral
$$\int_{(2)} L^{(2)}(u+1/2, \pi, \text{sym}^2) \frac{e^{u^2}}{u} \frac{G(u, \mu_{\pi})p^{(\frac{3}{2} + \alpha)u} }{
L_2(1/2-u, \tilde{\pi}, \text{sym}^2)}\frac{du}{2\pi i}.$$
Using Proposition \ref{prop1}, we can write this as a Dirichlet series. To this end, we write
$$L_2(v, \tilde{\pi}, \text{sym}^2)^{-1} = \sum_{j=0}^6 \frac{(-1)^j\overline{\lambda_{\pi}(2^j)}}{2^{jv}}$$
where (as one can check from \eqref{schur})
\begin{equation}\label{a}
\begin{split}
&\lambda_{\pi}(1) = 1, \quad \lambda(2) = A_{\pi}(1, 4), \quad \lambda_{\pi}(4) = A_{\pi}(2, 4), \quad \lambda_{\pi}(8) = A_{\pi}(8, 1) +A_{\pi}(1, 8), \\
& \lambda_{\pi}(16) = A_{\pi}(4, 2), \quad \lambda_{\pi}(32) = A_{\pi}(4, 1), \quad \lambda_{\pi}(64) = 1. 
\end{split}
\end{equation}

On the other hand, using Lemma \ref{polefree}, we can shift the contour to $\Re u = -1$ picking up at most a residue at $u = 0$. Then using the functional equation (recall Proposition \ref{prop2}) and opening the Dirichlet series we arrive at (recall \eqref{lstar})
\begin{displaymath}
\begin{split}
& L^{\ast}(1/2, \pi, \text{sym}^2) G(0, \mu_{\pi})  = \sum_{j=0}^6 \frac{(-1)^j \lambda_{\pi}(2^j)}{2^{j/2}} \sum_{\substack{a, b , c \text{ odd}\\ (a, p) = 1}} \frac{A_{\pi}(b^2, c^2)}{(a^3b^2c)^{1/2}}  V_{\pi}\Big(\frac{a^3b^2c}{2^jp^{3/2 + \alpha}}\Big)\\
&\quad\quad\quad\quad-\Big(p^{1/2} \overline{A_{\pi}(p, p)} - \frac{1}{p^{1/2}}\Big)\sum_{j=0}^6 \frac{(-1)^j \overline{\lambda_{\pi}(2^j)}}{2^{j/2}} \sum_{\substack{a, b , c \text{ odd}\\ (a, p) = 1}} \frac{ A_{\tilde{\pi}}(b^2, c^2)}{(a^3b^2c)^{1/2}}  W_{\pi}\Big(\frac{a^3b^2c}{2^jp^{3/2 - \alpha}}\Big)
\end{split}
\end{displaymath}
where  
\begin{displaymath}
\begin{split}
V_{\pi}(x) &= \int_{(2)} x^{-u} \frac{e^{u^2}}{u}  G(u, \mu_{\pi})\frac{du}{2\pi i}, \\
W_{\pi}(x) &= \int_{(2)} x^{-u} \frac{e^{u^2}}{u}  G(-u, \mu_{\pi}) \frac{L_{\infty}(u + 1/2, \tilde{\pi}, \text{sym}^2)}{L_{\infty}( 1/2-u, \pi, \text{sym}^2)}\frac{du}{2\pi i}, 
\end{split}
\end{displaymath}
Note that $G(-u, \mu_{\pi})$ cancels all poles of $L_{\infty}(u + 1/2, \tilde{\pi}, \text{sym}^2)$ in a wide strip. 
  We record the following properties of $V_{\pi}(x)$ and $W_{\pi}(x)$:
 \begin{itemize}
 \item as $x \rightarrow \infty$, the functions  $V_{\pi}$, $W_{\pi}$ are rapidly decaying,   more precisely, for every $A > 0$ there exists $B > 0$ such that $V_{\pi}(x) , W_{\pi}(x) \ll   \|  \mu_{\pi} \|^B x^{-A}$, uniformly in $\| \Re \mu_{\pi} \| \leq A_0$;
 \item as $x \rightarrow 0$, we have $V_{\pi}(x)  =G(0, \mu_{\pi})+ O(x \| \mu_{\pi} \|^{20})$ and $|W_{\pi}(x)|  =|G(0, \mu_{\pi}) |
 + O(x \| \mu_{\pi} \|^{20})$, uniformly in $\| \Re \mu_{\pi} \| \leq A_0$; 
 \item as a function of $\mu_{\pi}$ the  functions $V_{\pi}(x)$, $W_{\pi}(x)$ are Weyl group invariant, holomorphic in $\| \Re \mu_{\pi} \| \leq A_0$  and of polynomial growth in this strip.
   \end{itemize}

 
To display the following formulae in a compact fashion we write
\begin{equation}\label{compact}
 \sum_{j=0}^6 \frac{(-1)^j \lambda_{\pi}(2^j)}{2^{jv}} = \sum_{i, j, k} \eta_{ijk} \frac{A_{\pi}(2^i, 2^j)}{2^{kv}}
 \end{equation}
for suitable $\eta_{i j k} \in \{-1, 0, 1\}$.  
 For $(\m_1\m_2, p) = 1$ we conclude
 \begin{displaymath}
 \begin{split}
 A_{\pi}(\m_1, \m_2) &L^{\ast}(1/2, \pi, \text{sym}^2) G(0, \mu_{\pi})\\
 &= \sum_{i, j, k} \frac{\eta_{ijk}}{2^{k/2}} \Big[ \sum_{\substack{a, b , c \text{ odd}\\ (a, p) = 1}} \frac{ \overline{A_{\pi}(\m_2, \m_1)}A_{\pi}(2^ib^2, 2^jc^2)}{(a^3b^2c)^{1/2}}  V_{\pi}\Big(\frac{a^3b^2c}{2^kp^{3/2 + \alpha}}\Big) \\
& -p^{1/2} \sum_{\substack{a, b , c \text{ odd}\\ (a, p) = 1}} \frac{ A_{\pi}(\m_1, \m_2)\overline{A_{\pi}(p, p) A_{\pi}(2^ib^2, 2^jc^2)}}{(a^3b^2c)^{1/2}} W_{\pi}\Big(\frac{a^3b^2c}{2^kp^{3/2 - \alpha}}\Big)\\
& +   p^{-1/2} \sum_{\substack{a, b , c \text{ odd}\\ (a, p) = 1}} \frac{ A_{\pi}(\m_1, \m_2)\overline{A_{\pi}(2^ib^2, 2^jc^2)}}{(a^3b^2c)^{1/2}}   W_{\pi}\Big(\frac{a^3b^2c}{2^kp^{3/2 - \alpha}}\Big)\Big] .
  \end{split}
 \end{displaymath}

Before we continue, we observe that without loss of generality we may assume 
\begin{equation}\label{m}
   \m_1\m_2 \leq p^{\alpha/2} \leq p^{1/3},
   \end{equation}
say, otherwise there is nothing to prove. 
 
  We now modify the middle term in the penultimate long display  as follows. By the rapid decay of the weight function and the condition $\alpha < 1/2$, the contribution of $p\mid b$ and the contribution of $p^2 \mid c$ are negligible. 
   If $p\nmid bc$, then $
 A_{\pi}(p, p) A_{\pi}(2^ib^2, 2^jc^2) = A_{\pi}(2^ipb^2, 2^jpc^2)$. If $p\nmid b$, $c = pc'$ with $p\nmid c$, then
 $$  A_{\pi}(p, p) A_{\pi}(2^ib^2, 2^jc^2) =  A_{\pi}(p, p) A_{\pi}(1, p^2) A_{\pi}(2^ib^2, 2^j(c')^2) .$$
From Lemma \ref{ram} one verifies
 $$A_{\pi}(p, p)A_{\pi}(1, p^2) = A_{\pi}(p, p^3) + \frac{A_{\pi}(1, p^2)   + \overline{A_{\pi}(p^2, 1)}}{p}. $$
 Using bounds towards Ramanujan, we estimate trivially
 \begin{displaymath}
 \begin{split}
& p^{1/2} \sum_{\substack{a, b , c'\\ (abc', p) = 1}} \frac{ A_{\pi}( \m_1, \m_2)(\overline{A_{\pi}(1, p^2) }  +  A_{\pi}(p^2, 1))  \overline{  A_{\pi}(2^ib^2, 2^jc^2)}}{2^kp(a^3b^2pc')^{1/2}}  W_{\pi}\Big(\frac{a^3b^2c'}{p^{1/2 - \alpha}}\Big) \\
 &\preccurlyeq  \frac{(\m_1\m_2)^{ \theta_3 }}{p} p^{(\frac{1}{2} - \alpha) (\frac{1}{2}  + 2 \theta_3) }  \leq (\m_1\m_2)^{\theta_3} p^{ \theta_3 - \frac{3}{4} } 
 \end{split}
 \end{displaymath}
 where here and in the following all implied constants depend polynomially on $\mu_\pi$. 
 We conclude 
  \begin{equation}\label{beforeKuz}
 \begin{split}
  A_{\pi}(\m_1, \m_2) L^{\ast}(1/2, \pi, \text{sym}^2) G(0, \mu_{\pi}) = S_1(\pi) + S_2(\pi) + S_3(\pi) + O\big( (\m_1\m_2)^{\theta_3} p^{ \theta_3 - \frac{3}{4}+\varepsilon} \big) 
   \end{split}
 \end{equation}
 where
    \begin{equation}\label{sj}
 \begin{split}
S_1(\pi) = &\sum_{i, j, k} \eta_{ijk}  \sum_{\substack{a, b , c \text{ odd}\\ (a, p) = 1}} \frac{ \overline{A_{\pi}(\m_2, \m_1)}A_{\pi}(2^ib^2, 2^jc^2)}{(a^3b^2c)^{1/2}}  V_{\pi}\Big(\frac{a^3b^2c}{2^kp^{3/2 + \alpha}}\Big), \\
S_2(\pi) = & -p^{1/2} \sum_{i, j, k} \eta_{ijk} \sum_{\substack{a, b , c \text{ odd} \\ (a, p) = 1}} \frac{ A_{\pi}( \m_1, \m_2)\overline{ A_{\pi}(2^ipb^2, 2^jpc^2)}}{(a^3b^2c)^{1/2}}  W_{\pi}\Big(\frac{a^3b^2c}{2^kp^{3/2 - \alpha}}\Big),\\
S_3(\pi) = &   p^{-1/2} \sum_{i, j, k} \eta_{ijk} \sum_{\substack{a, b , c \text{ odd}\\ (a, p) = 1}} \frac{ A_{\pi}( \m_1, \m_2)\overline{A_{\pi}(2^ib^2, 2^jc^2)}}{(a^3b^2c)^{1/2}}   W_{\pi}\Big(\frac{a^3b^2c}{2^kp^{3/2 - \alpha}}\Big)   .
   \end{split}
 \end{equation}

\subsection{Adding the rest of the spectrum}

To enable an application of the Kuznetsov formula, we need to analyze the expressions $S_i(\pi)$, $i = 1, 2, 3$, when $\pi$ is replaced by an oldform or an Eisenstein series occuring in the spectral decomposition, in other words, if $A_{\pi}(m, n)$ is replaced with $B_{\pi}^{(j)}(m, n)$ for a level one representation $\pi$ or $B^{(j)}_{(s_1, s_2)}(m, n)$ for a minimal Eisenstein series or $A_{(\phi, s)}(m, n)$ or $B_{(\phi, s)}^{(j)}(m, n)$ for a maximal Eisenstein series. By Weyl's law, this part of the spectrum has smaller density, but is arithmetically somewhat subtle.  If one is interested only in upper bounds and has positivity, these contributions can be added trivially, but Theorem \ref{thm1} asks for an asymptotic formula and we do not have positivity, so we need to estimate the rest of the spectrum carefully. 

\medskip

We start with the hardest part, the estimation of $S_2(\pi)$ for $\pi$ a translate of a level 1 form $\pi$ as described in   Corollary \ref{cor4},  so we substitute $B_{\pi}^{(j)}(\ast, \ast)$ for $A_{\pi}(\ast, \ast)$ for $j = 0, 1, 2$. We call the corresponding quantity $S^{(j)}_2(\pi)$. 

By the rapid decay of $V_{\tilde{\pi}}$, the contribution $p \mid b$ is negligible, and also the contribution $p^2 \mid c$ is negligible. For $p^2\nmid c, p \nmid b$ we have by \eqref{hecke} that
\begin{equation}\label{usehecke}
\begin{split}
A_{\pi}&(2^ipb^2,2^jpc^2) = A_{\pi}(2^ipb^2, 2^jc^2) A_{\pi}(1, p) - A_{\pi}(2^ib^2, 2^jc^2) - \delta_{p \mid c}A_{\pi}(2^ip^2b^2, 2^jc^2/p)\\
& = A_{\pi}(2^i b^2, 2^jc^2)(A_{\pi}(p, 1) A_{\pi}(1, p) - 1) -   \delta_{p \mid c}A_{\pi}(2^ib^2, 2^j(c/p)^2)A_{\pi}(p^2, p)\\
& =  A_{\pi}(2^i b^2, 2^jc^2)A_{\pi}(p, p)
  -   \delta_{p \mid c}A_{\pi}(2^ib^2, 2^j(c/p)^2)A_{\pi}(p^2, p).
  \end{split}
\end{equation}
We conclude that 
\begin{equation*} 
\begin{split}
&S^{(0)}_2(\pi)  =  -p^{1/2} \sum_{i, j, k} \eta_{ijk} \sum_{ \substack{a, b, c\, \text{odd}\\ (a, p) = 1}} \frac{ B^{(0)}_{\pi}( \m_1, \m_2)\overline{ B^{(0)}_{\pi}(2^ipb^2, 2^jpc^2)}}{(a^3b^2c)^{1/2}}  W_{\pi}\Big(\frac{a^3b^2c}{2^kp^{3/2 - \alpha}}\Big)
\end{split}
\end{equation*}
equals up to a negligible error
\begin{equation}\label{compute} 
\begin{split}
&  -p^{1/2} A_{\pi}( \m_1, \m_2)\sum_{i, j, k} \eta_{ijk} \sum_{ \substack{a, b, c\, \text{odd}\\ (a, p) = 1}} \frac{\overline{A (2^ib^2, 2^jc^2)}}{(a^3b^2c)^{1/2}}\\
&\quad\quad\quad\quad\quad \quad\quad\quad\quad\quad\quad  \times  \Big[ \overline{ A_{\pi}(p, p)}W_{\pi}\Big(\frac{a^3b^2c}{2^kp^{3/2 - \alpha}}\Big) - \frac{\overline{A_{\pi}(p^2, p)}}{p^{1/2}}W_{\pi}\Big(\frac{a^3b^2c}{2^kp^{1/2 - \alpha}}\Big) \Big]\\
  =  & -p^{1/2} A_{\pi}(\m_1, \m_2)  \int_{(2)}  \frac{L^{(2)}(u + 1/2, \tilde{\pi}, \text{sym}^2)}{L_2(1/2 - u, \pi, \text{sym}^2)}  \frac{e^{u^2}}{u} G(-u, \mu_{\pi}) \frac{L_{\infty}(u + 1/2, \tilde{\pi}, \text{sym}^2)}{L_{\infty}( 1/2-u, \pi, \text{sym}^2)} \\&\quad\quad\quad\quad\quad \quad\quad\quad\quad\quad\quad \times \Big(\overline{A_{\pi}(p, p)} p^{(\frac{3}{2} - \alpha)u} -\frac{\overline{A_{\pi}(p^2, p)}}{p^{1/2}} p^{(\frac{1}{2} - \alpha)u} \Big) \frac{du}{2\pi i}.
\end{split}
\end{equation}
Shifting the contour to $\Re u = \varepsilon$ and recalling that $\pi$ has full level and hence its $L$-function is bounded independently in terms of $p$, we obtain
$$S^{(0)}_2(\pi)\preccurlyeq   p^{\frac{1}{2} + 2\theta_3}(\m_1\m_2)^{\theta_3}.  $$
For $p^2\nmid c, p \nmid b$ we have similarly
\begin{displaymath}
\begin{split}
A^{(1)}_{\pi}&(2^ipb^2,2^jpc^2)  = \sqrt{p} \big(A_{\pi}(2^ib^2,2^jp^2c^2)+ A_{\pi}(2^ipb^2,2^jc^2) \big)\\
& =  \sqrt{p}  \big( A_{\pi}(2^ib^2,2^jc^2)A_{\pi}(1, p^2) - \delta_{p \mid c}A_{\pi}(2^ipb^2,2^jc^2/p)A_{\pi}(1, p) + \delta_{p \mid c}A_{\pi}(2^ib^2,2^jc^2/p) \big)\\
&\quad  + \sqrt{p} \big(A_{\pi}(2^ib^2,2^jc^2)A_{\pi}(p, 1) - \delta_{p \mid c} A_{\pi}(2^ib^2,2^jc^2/p)\big)\\
& = \sqrt{p} \big( A_{\pi}(2^ib^2,2^jc^2)A_{\pi}(1, p)^2 - \delta_{p\mid c}  A_{\pi}(2^ib^2,2^j(c/p)^2)A_{\pi}(p, p) A_{\pi}(1, p)\big).
  \end{split}
\end{displaymath}
By the same argument as above we obtain
\begin{displaymath}
\begin{split}
S^{(1)}_2(\pi) & \ll  p^{1/2} \cdot |c_{10} A_{\pi}( \m_1, \m_2) | \\
& \quad\quad\quad \times \Big(|c_{10}|\Big( |A_{\pi}(p, p)| + \frac{|A_{\pi}(p^2, p)|}{p^{1/2}}\Big) + p^{1/2}\Big( |A_{\pi}(1, p)^2| + \frac{|A_{\pi}(p, p)A_{\pi}(1, p)|}{p^{1/2}}\Big)\Big) \\
& \preccurlyeq p^{\frac{1}{2} + \theta_3 - \frac{1}{2} + \frac{1}{2} + 2\theta_{3}  } (\m_1\m_2)^{\theta_3} = p^{\frac{1}{2} + 3\theta_3   } (\m_1\m_2)^{\theta_3}. 
\end{split}
\end{displaymath}
Finally, for  $p^2\nmid c, p \nmid b$ we have 
\begin{displaymath}
\begin{split}
A^{(2)}_{\pi}(2^ipb^2,2^jpc^2) & =pA_{\pi}(2^ib^2,2^jpc^2) = p \big(A_{\pi}(2^ib^2,2^jc^2)A_{\pi}(1, p) - \delta_{p\mid c} A_{\pi}(2^ipb^2,2^jc^2/p)\big)\\
& = p \big(A_{\pi}(2^ib^2,2^jc^2)A_{\pi}(1, p) - \delta_{p\mid c} A_{\pi}(2^ib^2,2^j(c/p)^2)A_{\pi}(p, p)\big)
   \end{split}
\end{displaymath}
and we obtain
\begin{displaymath}
\begin{split}
S^{(2)}_2(\pi) & \ll  p^{1/2} \cdot |c_{20} A_{\pi}(\m_1, \m_2) | \Big(|c_{20}|\Big( |A_{\pi}(p, p)| + \frac{|A_{\pi}(p^2, p)|}{p^{1/2}}\Big) \\
& \quad\quad + |c_{21}|p^{1/2}\Big( |A_{\pi}(1, p)^2| + \frac{|A_{\pi}(p, p)A_{\pi}(1, p)|}{p^{1/2}}\Big)    +p \Big(|A_{\pi}(1, p)|+ \frac{|A_{\pi}(p, p)|}{p^{1/2}}\Big)\Big) \\
& \preccurlyeq p^{1+ 2\theta_{3} } (\m_1 \m_2)^{\theta_3} .
\end{split}
\end{displaymath}
Combining the previous estimates, we obtain
\begin{equation}\label{error1}
\sum_{j=0}^2 S^{(j)}_2(\pi) \preccurlyeq p^{1+ 2\theta_{3} } (\m_1\m_2)^{\theta_3} .
\end{equation}

We can  run the same argument if $\pi$ is a minimal Eisenstein series or a maximal Eisenstein \emph{old}form, except that we can replace $\theta_3$ with $0$ or $\theta_2$ respectively. The key point is that even for Eisenstein series the integrand on the right hand side of  \eqref{compute} is holomorphic since the possible poles of $L^{(2)}(u+1/2, \tilde{\pi}, \text{sym}^2)$ at $$u \in \{ 1/2 + \overline{\mu_{\pi, j}} + \overline{\mu_{\pi, i}} \mid 1 \leq i, j \leq 3\}  = \{ 1/2 - \mu_{\pi, j} -  \mu_{\pi, i} \mid 1 \leq i, j \leq 3\} $$ are canceled by the denominator $L_{\infty}(1/2 - u, \pi, \text{sym}^2)$. Hence the same contour shift to $\Re u = \varepsilon$ is possible without collecting additional residues. 

\medskip

The estimation of  $S_i(\pi)$ for $i = 1,3$ and $\pi$ a cuspidal or Eisenstein oldform is similar, but much simpler and dominated by \eqref{error1}.  We obtain in a straightforward fashion
\begin{displaymath}
\begin{split}
& S_1^{(0)}(\pi) \ll  |A_{\pi}( \m_2, \m_1)|p^{\varepsilon} \preccurlyeq (\m_1\m_2)^{\theta_3}  ,\\
& S_1^{(1)}(\pi) \ll  |c_{10}|^2|A_{\pi}( \m_1, \m_2)| p^{\varepsilon} \preccurlyeq   p^{2\theta_3 - 1 } (\m_1\m_2)^{\theta_3 } ,\\
& S_1^{(2)}(\pi) \ll |c_{20}|^2|A_{\pi}( \m_1, \m_2)| p^{\varepsilon} \preccurlyeq   p^{2\theta_3 - 1 } (\m_1\m_2)^{\theta_3 } .
\end{split}
\end{displaymath}
The corresponding bounds for $S_3$ save an additional factor $p^{1/2}$. 

\medskip

We finally turn to maximal Eisenstein newforms, given by a pair $\pi = (\phi, s)$  as in Proposition \ref{prop5}a. Here we need two additional pieces of information. 

When we shift the contour in \eqref{compute} to $\Re u = \varepsilon$, we need to control the $L$-function $L^{(2)}(u, \tilde{\pi}, \text{sym}^2)$ by the convexity bound $O(p^{3/4+\varepsilon})$ on the line $\Re u = 1/2 + \varepsilon$. (To be able to apply the convexity bound, we need absolute convergence of the $L$-series to the right of $1$ which follows from \eqref{factor} in the cuspidal and non-cuspidal case.)  
For $S_2(\pi)$ we use the Hecke relations as in \eqref{usehecke} and 
observe from \eqref{satmax} that
\begin{displaymath}
\begin{split}
&A_{\pi}(p, p) = p^{-1} + p^{3s - \frac{1}{2}} \ll p^{-1/2}, 
\quad  A_{\pi}(p^2, p)  = p^{-3/2 + s} + p^{-1 + 4 s} \ll p^{-1}
\end{split}
\end{displaymath}
since $\phi$ is ramified at $p$, so $\{\alpha_p, \beta_p\} = \{p^{-1/2}, 0\}$. We then conclude by the same procedure as above  that
\begin{equation}\label{error2}
\begin{split}
|S_1(\pi)| + |S_2(\pi)|  + |S_3(\pi)| &\preccurlyeq |A_{\pi}(\m_1, \m_2)| p^{3/4 }.\\
\end{split}
\end{equation}

\subsection{Applying the Kuznetsov formula}

We choose $h \in \mathcal{H}(A_0)$ for some sufficiently large $A_0 > 10$ and let
\begin{equation}\label{tildeh}
\tilde{h}(\mu) = \frac{h(\mu)}{G(0, \mu)} \in \tilde{\mathcal{H}}(10).
\end{equation}
We define
\begin{equation}\label{ccusp}
\begin{split}
\mathfrak{C}_{\text{new}} 
&= c_{\text{cusp}} \sum_{\pi \in \mathcal{B}(p)}  \frac{S_1(\pi) + S_2(\pi) + S_3(\pi)}{L(1, \pi, \text{Ad})} \tilde{h}(\mu_{\pi}) . \end{split}
\end{equation}
By \eqref{beforeKuz} we have 
$$\mathfrak{C}_{\text{new}} =  c_{\text{cusp}} \sum_{\pi \in \mathcal{B}(p)} \frac{A_{\pi}(\m_1, \m_2) L^{\ast}(1/2, \pi, \text{sym}^2) h(\mu_{\pi})}{L(1, \pi, \text{Ad})}  + O\big( (\m_1\m_2)^{\theta_3} p^{2+ \theta_3 - \frac{3}{4}+\varepsilon}\big ) 
$$
where we used the bound
$$ \sum_{\pi \in \mathcal{B}(p)}  \frac{1}{L(1, \pi, \text{Ad})} \preccurlyeq p^{2}$$
for the error term which follows easily from the Kuznetsov formula (cf.\ e.g.\ \cite[Theorem 4]{BBM} for a much more refined result). We define in an analogous fashion the quantities
$$\mathfrak{C}_{\text{old}}, \quad \mathfrak{E}_{\text{min}}, \quad \mathfrak{E}_{\text{max, new}}, \quad \mathfrak{E}_{\text{max, old}}$$
as in \eqref{spec}. Using the lower bounds \cite{Br}
$$L(1, \pi,  \text{Ad}) \gg  \mu_{\pi}^{-A} $$
for some sufficiently large $A$ and a level 1 oldform $\pi$ on ${\rm PGL}(3)$, as well as the standard bounds
$$\zeta(1 + it) \gg_{\varepsilon} (1+ |t|)^{-\varepsilon}, \quad L(1 + i t, \phi) \gg ((1 + |t|)\text{cond}(\phi))^{-\varepsilon}, \quad L(1, \phi, \text{Ad}) \gg \text{cond}(\phi)^{-\varepsilon},$$
together with \eqref{error1} and \eqref{error2} and Weyl's law, we obtain
\begin{displaymath}
\begin{split}
& |\mathfrak{C}_{\text{old}}| + | \mathfrak{E}_{\text{min}}| +  |\mathfrak{E}_{\text{max, old}}| \preccurlyeq   p^{1+2\theta_3 } (\m_1\m_2)^{\theta_3},\\
&\mathfrak{E}_{\text{max, old}} \preccurlyeq p^{1 + \frac{3}{4} }  (\m_1\m_2)^{\theta_2}.
\end{split}
\end{displaymath}
We summarize the discussion by stating that
\begin{equation}\label{formula}
\begin{split}
& c_{\text{cusp}} \sum_{\pi \in \mathcal{B}(p)} \frac{A_{\pi}(\m_1, \m_2) L^{\ast}(1/2, \pi, \text{sym}^2) h(\mu_{\pi})}{L(1, \pi, \text{Ad})} \\
& = \mathfrak{C}_{\text{new}}  + \mathfrak{C}_{\text{old}}+ \mathfrak{E}_{\text{min}}+ \mathfrak{E}_{\text{max, new}}+ \mathfrak{E}_{\text{max, old}} + O\big(p^{7/4+\varepsilon} (\m_1\m_2)^{5/14}\big). 
\end{split}
\end{equation}
The main term on the right hand side is a form where the Kuznetsov formula  as in Lemma~\ref{lem.Kuznetsov}  can be applied. We call the corresponding terms on the geometric side $\mathfrak{D}, \mathfrak{S}_4, \mathfrak{S}_5, \mathfrak{S}_6$ and estimate each of them in turn.  According to the definition \eqref{ccusp}, each of these splits naturally into three summands, $\mathfrak{D}^{(j)}, \mathfrak{S}^{(j)}_4, \mathfrak{S}^{(j)}_5, \mathfrak{S}^{(j)}_6$, say, for $j = 1, 2, 3$, and we recall the definition of $S_j(\pi)$ in in \eqref{sj}. 

\medskip

We start with the three diagonal terms $\mathfrak{D}^{(j)}$. We have
$$\mathfrak{D}^{(1)} = P\sum_{i, j, k} \eta_{ijk}  \sum_{\substack{a, b , c \text{ odd}\\ (a, p) = 1}} \frac{ \delta_{\m_2 = 2^i b^2}\delta_{ \m_1 = 2^j c^2} }{(a^3b^2c)^{1/2}}  \int_{\Re \mu = 0} V_{\pi}\Big(\frac{a^3b^2c}{2^kp^{3/2 + \alpha}}\Big)  \tilde{h}(\mu) \text{spec}(\mu) d\mu.$$
Recalling \eqref{a} and \eqref{compact}, for odd $\m_1, \m_2$ this equals
\begin{displaymath}
\begin{split}
&\delta_{\substack{\m_1 = \square\\ \m_2 = \square}}\sum_{(a, p) = 1} \frac{ P }{a^{3/2} (\m_1\m_2^2)^{1/4}}  \int_{\Re \mu = 0} \Big(V_{\pi}\Big(\frac{a^3 \sqrt{\m_1\m_2^2}}{ p^{3/2 + \alpha}}\Big)  +  V_{\pi}\Big(\frac{a^3 \sqrt{\m_1\m_2^2}}{ 64p^{3/2 + \alpha}}\Big)\Big)\tilde{h}(\mu) \text{spec}(\mu) d\mu\\
= &\delta_{\substack{\m_1 = \square\\ \m_2 = \square}} \frac{P}{(\m_1\m_2^2)^{1/4}}  \int_{\Re \mu = 0} \int_{(2)} \frac{p^{(\frac{3}{2} + \alpha)u} (1 + 64^u)}{(\m_1^{1/2} \m_2)^{u}} \zeta^{(p)} (3u + \tfrac{3}{2}) \frac{e^{u^2}}{u} G(u, \mu_{\pi}) \frac{du}{2\pi i}  \tilde{h}(\mu) \text{spec}(\mu) d\mu.
\end{split}
\end{displaymath}
We shift the $u$-contour to $\Re u = -1/3$, say, picking up a pole  at $u = 0$. By \eqref{tildeh} we obtain
$$\delta_{\substack{\m_1 = \square\\ \m_2 = \square}}\frac{2 P\zeta  (  \tfrac{3}{2})(1 + O(p^{-1})) }{(\m_1\m_2^2)^{1/4}}  \int_{\Re \mu = 0}       h(\mu) \text{spec}(\mu) d\mu + O(Pp^{-1/2}) .$$

 We clearly have $\mathfrak{D}^{(2)} = 0$ since the delta-condition can never be satisfied. 
 
 The analysis of $\mathfrak{D}^{(3)}$ is almost identical to the  analysis of $\mathfrak{D}^{(1)}$, but saves a factor $p^{1/2}$, so that  $\mathfrak{D}^{(3)} \ll p^{3/2} $. All in all we conclude 
\begin{equation}\label{D}
\sum_{j=1}^3 \mathfrak{D}^{(j)} = \delta_{\substack{\m_1 = \square\\ \m_2 = \square}}\frac{Cp^2}{(\m_1\m_2^2)^{1/4}} + O(p^{3/2} )
\end{equation}
 for $C = 2 \zeta(3/2)\int_{\Re \mu = 0}       h(\mu) \text{spec}(\mu) d\mu$. 
 
 \subsection{Estimates for the $w_4$ and $w_5$ elements -- the easy cases}
 
We continue with bounds for $\mathfrak{S}^{(j)}_4$ and $\mathfrak{S}^{(j)}_5$. The only case that requires some work is $\mathfrak{S}^{(2)}_5$ which will be postponed to the next subsection. 
\medskip

Using the rapid decay of $V_{\pi}$ and  \eqref{lem8a} with $j=0$ and very large $A$, we have 
$$\mathfrak{S}^{(1)}_4 \ll  P \sum_{i, j \leq 3}   \sum_{\substack{a, b , c \text{ odd}\\ (a, p) = 1}} \frac{1}{(a^3b^2c)^{1/2}} \sum_{\pm} \sum_{\substack{p D_2 \mid D_1\\ \m_1 D_1 = 2^i b^2 D_2^2\\ a^3 b^2 c \preccurlyeq p^{3/2 + \alpha}\\ D_1D_2 \preccurlyeq \m_1\m_2c^2 }} \frac{|\tilde{S}(\mp 2^j c^2, \m_1, \m_2; D_2, D_1)|}{D_1D_2}  + O(p^{-100}).$$
 By Lemma \ref{lem1}(d) we obtain
$$\mathfrak{S}^{(1)}_4 \preccurlyeq P \sum_{i  \leq 3}   \sum_{\substack{a, b , c \text{ odd}\\ (a, p) = 1}} \frac{1}{(a^3b^2c)^{1/2}}   \sum_{\substack{p D_2 \mid D_1\\ \m_1 D_1 = 2^i b^2 D_2^2\\ a^3 b^2 c \preccurlyeq p^{3/2 + \alpha}\\ D_1D_2 \preccurlyeq \m_1\m_2c^2 }} \frac{\m_2D_2}{D_1} + O(p^{-100}).$$
We write $pD_2 \delta = D_1$ so that  $  \delta = 2^i b^2 D_2/(\m_1 p)$, and hence in particular $p \mid b^2 D_2$. The size conditions exclude the  case $p \mid b$ for $\alpha < 1/2$, so that we obtain
$$ \mathfrak{S}^{(1)}_4 \preccurlyeq  p^{2}     \sum_{\substack{a^3b^2 c \preccurlyeq p^{3/2 + \alpha }\\  b^2 D_2^3 \preccurlyeq \m_1^2\m_2 c^2 p \\ p \mid D_2} } \frac{1}{(a^3b^2c)^{1/2}}    \frac{\m_1\m_2}{b^2 D_2}+ O(p^{-100}).$$
Writing $D_2 = p \delta'$, it is easy to see that 
$$ \mathfrak{S}^{(1)}_4 \preccurlyeq\m_1\m_2 p^{2 } \cdot \frac{1}{p} \cdot p^{\frac{1}{2}(\frac{3}{2} + \alpha)} =\m_1\m_2 p^{\frac{7}{4} + \frac{\alpha}{2} }. $$

The bounds $\mathfrak{S}^{(2)}_4$ and $\mathfrak{S}^{(3)}_4$ are easier, since the summation conditions force the contribution to be negligible. For $\mathfrak{S}^{(2)}_4$ we are summing, up a negligible error, over quintuples $(a, b, c, D_1, D_2)$ satisfying
$$p D_2 \mid D_1, \quad p2^j c^2 D_1 = \m_1 D_2^2, \quad  a^3 b^2 c \preccurlyeq p^{3/2 - \alpha}, \quad  D_1D_2 \preccurlyeq \m_2b^2c^2 p^{2}.$$
With  $pD_2 \delta = D_1$  we have $p^2 2^j c^2 \delta = \m_1 D_2$, so  that the last condition becomes $c^2 \delta^3 p^3 \preccurlyeq \m_1^2\m_2b^2  \preccurlyeq \m_1^2\m_2 p^{3/2 - \alpha }$, which is impossible for $\m_1\m_2 \leq p^{1/3}$. 

Similarly, for $\mathfrak{S}^{(3)}_4$ we are summing, up a negligible error, over quintuples $(a, b, c, D_1, D_2)$ satisfying
$$p D_2 \mid D_1, \quad   2^j c^2 D_1 = \m_1 D_2^2, \quad  a^3 b^2 c \preccurlyeq p^{3/2 - \alpha}, \quad  D_1D_2 \preccurlyeq\m_2b^2c^2$$
With  $pD_2 \delta = D_1$  we have $p 2^j c^2 \delta = \m_1D_2$, so  that the last condition becomes $c^2 \delta^3 p^3 \preccurlyeq \m_1^2\m_2 b^2  \preccurlyeq \m_1^2\m_2 p^{3/2 - \alpha}$,  which is impossible for $\m_1\m_2 \leq p^{1/3}$. 

We conclude
\begin{equation}\label{s4}
\sum_{j=1}^3 |\mathfrak{S}^{(j)}_4| \preccurlyeq \m_1\m_2 p^{\frac{7}{4} + \frac{\alpha}{2} }
\end{equation}
provided that \eqref{m} holds. 

\medskip

We now turn to the estimation for  $\mathfrak{S}^{(j)}_5$. For $j=1$, up to a negligible error, we are summing over quintuples $(a, b, c, D_1, D_2)$ satisfying
$$p \mid D_1 \mid D_2, \quad  \m_2D_2 = 2^j c^2 D_1^2, \quad  a^3 b^2 c \preccurlyeq p^{3/2 + \alpha}, \quad  D_1D_2 \preccurlyeq \m_1\m_2 c^2 .$$
We write $D_1 = pd_1$, $D_2 = pd_1\delta$, so that $\delta = 2^j c^2 p d_1/\m_2$, and the last condition  becomes $d_1^3 p^3 \preccurlyeq\m_1\m_2^2 $, which is impossible. 

We argue similarly for $j=3$. Here we have the conditions
$$p \mid D_1 \mid D_2, \quad  2^ib^2 D_2 =  \m_2D_1^2, \quad  a^3 b^2 c \preccurlyeq p^{3/2- \alpha}, \quad  D_1D_2 \preccurlyeq\m_1  b^2c^2 .$$
Again we write   $D_1 = pd_1$, $D_2 = pd_1\delta$, so that $pd_1 = 2^i b^2 \delta/\m_2$. The first size condition makes $p \mid b^2$ impossible, so that  $p \mid \delta$, say $\delta = p \delta'$ and $d_1 = 2^i b^2 \delta'/\m_2$, and the last size condition becomes $b^2 (\delta')^3p^3  \preccurlyeq \m_1\m_2^2 c^2 $. Since $\m_1\m_2 \leq p^{\alpha/2}$ with $\alpha > 0$, this is impossible. 

The case $j=2$ takes a little longer and is the subject of the next subsection.

\subsection{Estimation of $\mathfrak{S}^{(2)}_5$}

Arguing as before, we have, up to a negligible error,
\begin{displaymath}
\begin{split}
\mathfrak{S}^{(2)}_5 \ll  Pp^{1/2} \sum_{i, j \leq 3}\sum_{\pm}  \Big|   \sum_{\substack{a, b , c \text{ odd}\\ (a, p) = 1}} \sum_{\substack{p \mid D_1 \mid D_2\\ 2^ipb^2 D_2 =  \m_2D_1^2\\ a^3 b^2 c \preccurlyeq p^{3/2- \alpha} }} \frac{\tilde{S}(\mp \m_1, 2^ipb^2, 2^jpc^2; D_1, D_2)}{(a^3b^2c)^{1/2}D_1D_2}  \Phi_{w_5}\Big(\frac{\pm  2^{i+j} \m_1b^2 c^2 p^2}{D_1D_2} \Big) \Big|. 
\end{split}
\end{displaymath}
Again we write   $D_1 = pd_1$, $D_2 = pd_1\delta$, so that $d_1 = 2^i b^2 \delta/\m_2$. 
With this parametrization we obtain
\begin{equation}\label{s52}
\begin{split}
\mathfrak{S}^{(2)}_5 \ll  &Pp^{1/2} \sum_{i, j \leq 3}\sum_{\pm}  \Big|   \sum_{\substack{a, b , c \text{ odd}\\ (a, p) = 1\\   a^3 b^2 c \preccurlyeq p^{3/2- \alpha}}} \\
&\times \sum_{\delta \equiv 0 \, (\text{mod } \frac{\m_2}{(\m_2, b^2)})} \frac{\tilde{S}(\mp \m_1, 2^ipb^2, 2^jpc^2; \frac{2^ipb^2\delta}{\m_2}, \frac{2^ipb^2\delta^2}{\m_2})}{a^{3/2} b^5 c^{1/2}p^2   \delta^3 \m_2^{-2}}  \Phi_{w_5}\Big(\frac{\pm  2^{j-i} \m_1\m_2^2 c^2  }{  \delta^3 b^2} \Big) \Big|. 
\end{split}
\end{equation}
By definition, we have
\begin{equation}\label{ev-klo}
\begin{split}
\tilde{S}\Big(\mp &\m_1, 2^ipb^2, 2^jpc^2; \frac{2^ipb^2\delta}{\m_2}, \frac{2^ipb^2\delta^2}{\m_2}\Big)\\
& = \sum_{\substack{C_1\, ( \text{mod } 2^ipb^2\delta/\m_2)\\ (C_1, 2^ipb^2\delta/\m_2 ) = 1}} \sum_{ \substack{C_2 \, (\text{mod } 2^ipb^2\delta^2/\m_2) \\(C_2, \delta )=1}} e\left(  \m_2 \frac{\bar{C}_1 C_2}{ \delta}+ 2^jpc^2 \frac{\bar{C}_2}{\delta}\mp \m_1\m_2 \frac{C_1}{2^ipb^2\delta}\right)\\
& = \frac{\phi(2^ipb^2\delta/\m_2)}{\phi(2^ipb^2\delta)} \frac{\phi(2^i p b^2 \delta^2/\m_2)}{\phi(\delta)}
\\
 & \quad\quad\quad  \times \sum_{\substack{C_1\, ( \text{mod } 2^ipb^2\delta)\\ (C_1, 2^ipb^2\delta ) = 1}} \sum_{ \substack{C_2 \, (\text{mod } \delta) \\(C_2, \delta )=1}}e\left(  \m_2 \frac{\bar{C}_1 C_2}{ \delta}+ 2^jpc^2 \frac{\bar{C}_2}{\delta}\mp \m_1\m_2 \frac{C_1}{2^ipb^2\delta}\right). 
\end{split}
\end{equation}
The bound \eqref{lem8a} with $j=0$ and large $A$ implies, up to a negligible error,
\begin{equation}\label{negl}
\delta^3 b^2\preccurlyeq \m_1\m_2 c^2 \preccurlyeq p^{3 - 2\alpha } p^{\alpha/2}
\end{equation}
so that $p \nmid \delta$. 
Since also $p \nmid b$ and $(p, 2\m_2) = 1$, the $C_1$-sum splits off a Ramanujan sum modulo $p$, so the Kloosterman sum equals
\begin{equation}\label{klooster}
\begin{split}
&-\frac{\phi(2^ipb^2\delta/\m_2)}{\phi(2^ipb^2\delta)} \frac{\phi(2^i p b^2 \delta^2/\m_2)}{\phi(\delta)}\\
&\times 
 \sum_{\substack{C_1\, ( \text{mod } 2^i b^2\delta)\\ (C_1, 2^i b^2\delta ) = 1}} \sum_{ \substack{C_2 \, (\text{mod } \delta) \\(C_2, \delta )=1}}e\left( \m_2  \frac{\bar{C}_1 C_2}{ \delta}+ 2^jpc^2 \frac{\bar{C}_2}{\delta}\mp \m_1\m_2 \frac{C_1\bar{p}}{2^ib^2\delta}\right).
\end{split}
\end{equation}
We plug this back into \eqref{s52}. 

In order to prepare for Poisson summation in $c$, we use a smooth partition of unity,      split the $c$-sum into dyadic ranges $c \asymp C$ using a weight function $W$ and call the corresponding pieces $\mathfrak{S}_5^{(2)}(C)$. By  Poisson summation we have
\begin{displaymath}
\begin{split}
&\sum_{c \text{ odd}} W\Big( \frac{c}{C}\Big) \frac{1}{c^{1/2}}   \Phi_{w_5}\Big(\frac{\pm  2^{j-i} \m_1\m_2^2  c^2  }{  \delta^3 b^2} \Big)  e\Big( \frac{2^j pc^2 \bar{C}_2}{\delta}\Big)\\
 = & \sum_{\kappa = 0, 1} (-1)^{\kappa} \sum_{\gamma \, (\text{mod } \delta)} e\Big( \frac{2^{j+2\kappa} p\gamma^2 \bar{C}_2}{\delta}\Big) \frac{1}{\delta} \sum_{ c} e\Big(\frac{c\gamma}{\delta}\Big) \\
& \quad\quad\quad\quad\quad \times \int W\Big( \frac{2^{\kappa} x}{C}\Big) \frac{1}{2^{\kappa/2}x^{1/2}}   \Phi_{w_5}\Big(\frac{\pm  2^{j-i+2\kappa} \m_1\m_2^2  x^2  }{  \delta^3 b^2} \Big)  e\Big( - \frac{xc}{\delta}\Big) dx.
\end{split}
\end{displaymath}
Repeated integration by parts with \eqref{lem8a} shows that the $x$-integral is negligible unless
$$c \preccurlyeq  C^{\ast} :=  \frac{\delta}{C} + \frac{ (\m_1\m_2^2)^{1/3}}{C^{1/3}b^{2/3}} ,$$
otherwise we estimate the integral trivially by 
$$\preccurlyeq C^{1/2} \Big(\frac{\m_1\m_2^2 C^2}{\delta^3 b^2}\Big)^{-1/6}= \frac{\delta^{1/2} C^{1/6} b^{1/3}}{(\m_1\m_2^2)^{1/6}}$$ 
again using \eqref{lem8a} with $j=0$. Recalling \eqref{ev-klo} and \eqref{negl}, we obtain altogether
\begin{equation}\label{s25}
\mathfrak{S}^{(2)}_5(C) \preccurlyeq p^{5/2}\sum_{i, j  \leq 3} \sum_{\pm}  \sum_{\substack{a, b, \delta\\ a^3 b^2 C \preccurlyeq p^{3/2 - \alpha }\\ \delta^3 b^2 \preccurlyeq \m_1\m_2^2 C^2}} \frac{1}{a^{3/2} b^5 p^2\delta^3\m_2^{-2}} \sum_{c \ll C^{\ast}} \frac{\delta^{1/2} C^{1/6} b^{1/3}}{(\m_1\m_2^2)^{1/6}} 
\frac{1}{\m_2} \frac{pb^2\delta}{\m_2}\frac{1}{\delta} |\mathfrak{C}(\delta)|
\end{equation}
where
$$\mathfrak{C}(\delta) =\sum_{\substack{C_1\, ( \text{mod } 2^i b^2\delta)\\ (C_1, 2^i b^2\delta ) = 1}} \sum_{ \substack{C_2 \, (\text{mod } \delta) \\(C_2, \delta )=1}}\sum_{\gamma \, (\text{mod } \delta)} e\left(  \m_2 \frac{\bar{C}_1 C_2}{ \delta}+ 2^jp\gamma^2 \frac{\bar{C}_2}{\delta}\mp \m_1\m_2 \frac{C_1\bar{p}}{2^ib^2\delta} + \frac{c\gamma}{\delta}\right).$$
We need essentially best-possible bounds for this multiple character sum. If $(b, 2\delta \m_1\m_2) = 1$, we could factor out a Ramanujan modulo $b^2$ which is bounded, and then hope for square-root cancellation in all three  sums modulo $\delta$. The next lemma establishes this up to some common divisors. 
\begin{lemma} We have
$$\mathfrak{C}(\delta) \ll (b\delta)^{\varepsilon} \delta^{3/2} ((b^2 \delta, \m_1\m_2)(\delta^2 b^2, \delta \m_1\m_2, \m_1\m_2)(\delta, \m_2^{\infty}))^{1/2}.$$ 
\end{lemma}

\textbf{Proof.} We can reduce the triple sum  in an elementary way to a one-dimensional exponential sum as follows. 
Changing variables $\gamma \mapsto \gamma C_2$ and evaluating a Ramanujan sum, we have
\begin{displaymath}
\begin{split}
\mathfrak{C}(\delta) &= 
\sum_{\substack{C_1\, ( \text{mod } 2^i b^2\delta)\\ (C_1, 2^i b^2\delta ) = 1}} \sum_{ \substack{C_2 \, (\text{mod } \delta) \\(C_2, \delta )=1}}\sum_{\gamma \, (\text{mod } \delta)} e\left(  \m_2 \frac{\bar{C}_1 C_2}{ \delta}+ 2^jp\gamma^2 \frac{ C_2}{\delta}\mp \m_1\m_2 \frac{C_1\bar{p}}{2^ib^2\delta} + \frac{c\gamma C_2}{\delta}\right)\\
& = 
\sum_{\substack{C_1\, ( \text{mod } 2^i b^2\delta)\\ (C_1, 2^i b^2\delta) = 1}} \sum_{\gamma \, (\text{mod } \delta)} e\left(   \mp \m_1\m_2 \frac{C_1\bar{p}}{2^ib^2\delta} \right)r_{\delta}(\m_2 + 2^j p \gamma^2 C_1+ c \gamma C_1)\\
& = 
\sum_{d \mid \delta} d \mu\Big( \frac{\delta}{d}\Big) \underset{\m_2 + 2^j p \gamma^2 C_1+ c \gamma C_1 \equiv 0 \, (\text{mod }d)}{ \sum_{\substack{C_1\, ( \text{mod } 2^i b^2\delta)\\ (C_1, 2^i b^2\delta ) = 1}} \sum_{\gamma \, (\text{mod } \delta)} }e\left(   \mp \m_1\m_2 \frac{C_1\bar{p}}{2^ib^2\delta} \right).
\end{split}
\end{displaymath}
Since $(C_1, d) = 1$, we must have $(\m_2, d) = (2^jp\gamma^2 + c\gamma, d) =: g$, say, and the congruence becomes
$$C_1  \equiv - \frac{\m_2}{g} \overline{\frac{2^jp\gamma^2 + c\gamma}{g}} \, \Big(\text{mod } \frac{d}{g}\Big).$$
An exercise in M\"obius inversion allows us  to evaluate the $C_1$-sum getting
\begin{displaymath}
\begin{split}
\mathfrak{C}(\delta) &=  \sum_{ d \mid \delta} d \mu\Big( \frac{\delta}{d}\Big)\sum_{\substack{\gamma \, (\text{mod } \delta)\\ (2^j p \gamma^2 + c\gamma, d) = g}} \sum_{\substack{f_1f_2 =   2^i b^2 \delta\\ (f_1, d/g) = 1\\ f_2 \mid \m_1\m_2}} \mu(f_1) f_2 e\Big( \pm  \frac{\m_1\m_2}{f_2} \frac{\overline{p\frac{2^jp\gamma^2 + c\gamma}{g}f_1} \frac{\m_2}{g}}{d/g} \Big)
\end{split}
\end{displaymath}
with $g = (\m_2, d)$. We change variables $\gamma \mapsto \gamma \bar{p}$ and observe that the sum over $\gamma$ depends only on $\gamma$ modulo $d$. Thus we obtain
$$\mathfrak{C}(\delta) \ll \delta  \sum_{d \mid \delta} \sum_{\substack{f_1f_2 =   2^i b^2 \delta\\ (f_1, d/g) = 1\\ f_2 \mid \m_1\m_2}} f_2 \Big|\sum_{\substack{\gamma \, (\text{mod } d)\\ (2^j  \gamma^2 + c\gamma, d) = g}} e\Big( \pm  \frac{\m_1\m_2}{f_2} \frac{\overline{\frac{2^j\gamma^2 + c\gamma}{g}f_1} \frac{\m_2}{g}}{d/g} \Big)\Big|.$$
We decompose $d = d_1d_2$ where $(d_1, \m_2) = 1$ and $d_2 \mid \m_2^{\infty}$, and accordingly split the $\gamma$-sum. We estimate the $d_2$-part trivially. The $d_1$-part is given by
\begin{displaymath}
\begin{split}
& \sum_{\substack{\gamma \, (\text{mod } d_1)\\ (2^j  \gamma^2 + c\gamma, d_1) = 1}} e\Big( \pm  \frac{\m_1\m_2}{f_2} \frac{\overline{\bar{g}(2^j\gamma^2 + c\gamma)f_1} \frac{\m_2}{g} \overline{\frac{d_2}{g}}}{d_1} \Big)\\
& \ll \Big( d_1, \frac{\m_1\m_2}{f_2}\Big)\Big|  \sum_{\substack{\gamma \, (\text{mod } \frac{d_1}{(d_1, \m_1\m_2/f_2))}\\ (2^j  \gamma^2 + c\gamma, d_1) = 1}} e\Big( \pm  \frac{\m_1\m_2/f_2}{(d_1, \m_1\m_2/f_2)} \frac{\overline{\bar{g}(2^j\gamma^2 + c\gamma)f_1} \frac{\m_2}{g} \overline{\frac{d_2}{g}}}{d_1/(d_1, \m_1\m_2/f_2)} \Big)\Big|.
\end{split}
\end{displaymath}

The $\gamma$-sum is now a one-dimensional exponential sum with a rational function. We can split the modulus into prime powers $p^{\alpha}$. If $\alpha = 1$, we use the Riemann hypothesis over finite fields (cf.\ e.g. \cite{Bo} with $n=1$), while for $\alpha > 1$ the argument is elementary \cite[Section 12]{IK}. In either case, the $\gamma$-sum is bounded $(d_1/(d_1, \m_1\m_2/f_2))^{1/2 + \varepsilon}$, and we obtain
$$\mathfrak{C}(d) \ll \delta \sum_{d\mid \delta} \sum_{f_2 \mid (2^i b^2 \delta, \m_1\m_2)} f_2 \Big(d, \frac{\m_1\m_2}{f_2}\Big)\Big(\frac{d}{(d, \m_1\m_2/f_2)}\Big)^{1/2 + \varepsilon} (d, \m_2^{\infty})^{1/2},$$
which implies the desired bound.\\

We substitute the character sum bound into \eqref{s25} and obtain
\begin{displaymath}
\begin{split}
\mathfrak{S}^{(2)}_5(C) & \preccurlyeq p^{3/2}   \sum_{\substack{a, b, \delta\\ a^3 b^2 C \preccurlyeq p^{3/2 - \alpha}\\ \delta^3 b^2 \preccurlyeq \m_1\m_2 C^2 }} \frac{1}{a^{3/2} b^3 \delta^3} \sum_{c\preccurlyeq C^{\ast}} \frac{\delta^{1/2} C^{1/6} b^{1/3}}{(\m_1\m_2^2)^{1/6}} \delta^{3/2} \m_1\m_2 (\delta, \m_2^{\infty})^{1/2}\\
&\preccurlyeq p^{3/2}   \sum_{\substack{  b, \delta\\  b^2 C \preccurlyeq p^{3/2 - \alpha }\\ \delta^3 b^2 \preccurlyeq \m_1\m_2 C^2 }} \frac{(\m_1\m_2)^{5/6}}{  b^{8/3} \delta }  C^{1/6} \Big(1 + \frac{\delta}{C} + \frac{ (\m_1\m_2^2)^{2/3}}{C^{1/3}b^{2/3}}\Big)(\delta, \m_2^{\infty})^{1/2} \\
&\preccurlyeq p^{3/2 } C^{1/6} (\m_1\m_2)^{5/6}\Big( 1 + \frac{(\m_1\m_2^2)^{1/3}}{C^{1/3}}\Big).
\end{split}
\end{displaymath}
For $C \preccurlyeq p^{3/2 - \alpha}$ this gives the final bound
\begin{equation}\label{S5}
\mathfrak{S}^{(2)}_5(C)  \preccurlyeq p^{7/4 } (\m_1\m_2)^{5/2}.
\end{equation}

 \subsection{Estimates for the long Weyl element} Here we bound $\mathfrak{S}^{(j)}_6$ for $j = 1, 2, 3$. The hardest case is $j=1$, so let us first quickly dispense with the cases $j=2, 3$.

       \medskip
    
Up to a negligible error,    we have
\begin{displaymath}
\begin{split}
\mathfrak{S}^{(2)}_6 \ll Pp^{1/2} \sum_{i, j \leq 3}  \sum_{\epsilon_1, \epsilon_2 = \pm 1} \Big|  \sum_{\substack{a, b , c \text{ odd}\\ (a, p) = 1\\ a^3 b^2 c \preccurlyeq p^{3/2 - \alpha}}} & \sum_{  p \mid D_1, p \mid D_2 } \frac{S^{(p)}(\epsilon_2 \m_2, \epsilon_1\m_1, 2^i pb^2, 2^j pc^2; D_1, D_2)|}{(a^3b^2c)^{1/2}D_1D_2}\\
&\times  \Phi_{w_6}\Big(- \frac{\epsilon_2 \m_22^ipb^2D_2}{D_1^2}, - \frac{\epsilon_1 2^j pc^2 \m_1D_1}{D_2^2} \Big)\Big|. 
\end{split}
\end{displaymath}
  The bound \eqref{lem8b} with $j_1 = j_2 = 0$ and very large $A$ implies that we can truncate the $D_1, D_2$ sum, up to a negligible error, at
       $$ D_1 \preccurlyeq b^{4/3}(c\m_2)^{2/3}\m_1^{1/3}   p, \quad  D_2 \preccurlyeq c^{4/3} (b\m_1)^{2/3}\m_2^{1/3} p.$$
 We write $D_1 = p \delta_1$, $D_2 = p\delta_2$, so that by \eqref{m} we have $\delta_1 \preccurlyeq p^{1 - \alpha/3 }$ and $\delta_2       \preccurlyeq p^{2 - 2\alpha/3 }$. In particular, $p \nmid \delta_1$ and $p^2 \nmid \delta_2$. Let us now consider the Kloosterman sum
 $$S^{(p)}(\epsilon_2 \m_2, \epsilon_1\m_1, 2^i pb^2, 2^j pc^2; p \delta_1, p \delta_2).$$
    Let us first assume $p \nmid \delta_2$. Then by Lemma \ref{lem1}(a) it splits off a factor
    $$S^{(p)}(\epsilon_2 \bar{\delta}_1^2 \delta_2 \m_2, \epsilon_1\m_1 \bar{\delta}_2^2\delta_1, 2^i pb^2, 2^j pc^2, p, p) = 0$$
 by Lemma \ref{lem1}(b). Let us next assume that $p \parallel \delta_2$, say $\delta_2 = p \delta_2'$ with $p \nmid \delta_2$.  Then again by   Lemma \ref{lem1}(a) the Kloosterman sum splits off a factor
 $$S^{(p)}(\epsilon_2 (\bar{\delta}_1^2 \delta'_2 \m_2, \epsilon_1\m_1 (\bar{\delta}'_2)^2\delta_1, 2^i pb^2, 2^j pc^2, p, p^2) $$
 which   by Lemma \ref{klooster-lemma}(a) splits off a factor $S ( \epsilon_1\m_1 (\bar{\delta}'_2)^2\delta_1, 2^j p^2 c^2, p^2) = 0$ since $p \nmid p \delta_1$. 
 
 \medskip
 
Next, for $\mathfrak{S}^{(3)}_6$ we get the size conditions
 $$a^3 b^2 c \preccurlyeq p^{3/2 - \alpha}, \quad    D_1 \preccurlyeq b^{4/3}(c\m_2)^{2/3}\m_1^{1/3}   , \quad D_2 \preccurlyeq  c^{4/3} (b\m_1)^{2/3}\m_2^{1/3},$$
so in particular
$$p \mid D_1 \preccurlyeq  (\m_2^2\m_1)^{1/3}  p^{1 - 2\alpha/3 } \ll p^{1 - \alpha/3 }$$
which is impossible. 
 

\medskip

Finally, up to a negligible error we have 
\begin{displaymath}
\begin{split}
\mathfrak{S}^{(1)}_6 \ll P \sum_{i, j \leq 3}  \sum_{\epsilon_1, \epsilon_2 = \pm 1} \Big|  \sum_{\substack{a, b , c \text{ odd}\\ (a, p) = 1\\ a^3 b^2 c \preccurlyeq p^{3/2 + \alpha}}} & \sum_{  p \mid D_1, p \mid D_2 } \frac{S^{(p)}(\epsilon_2 2^j c^2, \epsilon_12^i b^2, \m_2, \m_1; D_1, D_2)}{(a^3b^2c)^{1/2}D_1D_2}\\
& \times \Phi_{w_6}\Big(- \frac{\epsilon_2 2^jc^2\m_2D_2}{D_1^2}, - \frac{\epsilon_1 2^i b^2 \m_1D_1}{D_2^2} \Big)\Big|. 
\end{split}
\end{displaymath}
An accurate analysis of this term takes some time and needs again an application of Poisson summation. The final bound will be given in \eqref{S6} below.

As before, we use the parametrization $D_1 = p \delta_1$, $D_2 = p \delta_2$. 
The bound \eqref{lem8b} with $j_1 = j_2 = 0$ and very large $A$ implies that we can truncate the $\delta_1, 
 \delta_2$ sum, up to a negligible error, at
\begin{equation}\label{bound-delta}
 p\delta_1 \preccurlyeq c^{4/3}(b\m_2)^{2/3}\m_1^{1/3}  , \quad p \delta_2 \preccurlyeq b^{4/3} (\m_1c)^{2/3}\m_2^{1/3}.
 \end{equation}
  The size restrictions imply in particular $p \nmid \delta_2$ and $p^2 \nmid \delta_1$. 
  
  Let us first assume $p \mid \delta_1$, say $\delta_1 = p \delta_1'$ with $p \nmid \delta_1'$. Then by Lemma \ref{lem1}(a) and Lemma \ref{klooster-lemma}(b)  the Kloosterman sum  splits off a factor $$S^{(p)}(\epsilon_2 2^j c^2 (\bar{\delta}_1')^2\delta_2, \epsilon_12^i b^2\bar{\delta}_2^2 \delta_1', \m_2, \m_1; p^2 , p ) = - S(\epsilon_2 2^j c^2 (\bar{\delta}_1')^2\delta_2, \m_2p, p^2) = p \delta_{p \mid c}.$$
We write $c = pc'$ and estimate the remaining piece of the Kloosterman sum trivially using the trivial bound $S^{(p)}(*, *, *, *, D_1, D_2) \ll (D_1D_2)^{1+\varepsilon}$. 
Thus in this subcase we are left with bounding 
  \begin{displaymath}
\begin{split}
 &   p^{2}    \sum_{   b^2 c' \preccurlyeq p^{1/2 + \alpha}}  \sum_{  \delta_1' \preccurlyeq (\frac{(c')^{4}b^2\m^2_2\m_1}{p^2})^{1/3} } \sum_{\delta_2 \preccurlyeq(\frac{b^{4} \m_1^2(c')^{2}\m_2}{p})^{1/3} }\frac{1}{ b(c')^{1/2} p^{5/2}}\\
& \preccurlyeq p^{2}  \sum_{   b^2 c \preccurlyeq p^{1/2 + \alpha}} \frac{1}{ bc^{1/2} p^{5/2}} \frac{c^2 b^2 \m_1\m_2}{p} \preccurlyeq\m_1\m_2  p^{-1/4 + 5\alpha/2 },
\end{split}
\end{displaymath}
which is easily dominated by previous error terms. 

Let us now restrict to the case $p \nmid \delta_1$. Then by Lemma \ref{lem1}(a \& b), together with the fact that we can assume $p \nmid b$, we are left with bounding
 \begin{displaymath}
\begin{split}
p    \sum_{   b^2 c \preccurlyeq p^{3/2 + \alpha}  } & \sum_{  p \nmid \delta_1\delta_2 } \frac{S^{(1)}(\epsilon_2 2^j c^2\bar{p}, \epsilon_12^i b^2\bar{p}, \m_2, \m_1; \delta_1, \delta_2)}{ bc^{1/2}\delta_1\delta_2} \Phi_{w_6}\Big(- \frac{\epsilon_2 2^jc^2\m_2\delta_2}{p\delta_1^2}, - \frac{\epsilon_1 2^i b^2 \m_1\delta_1}{p\delta_2^2} \Big). 
\end{split}
\end{displaymath}
We apply a smooth partition of unity to restrict into dyadic ranges $c \asymp C$  and  apply the decomposition Lemma \ref{lem1}(c) getting expressions of the form 
 \begin{displaymath}
\begin{split}
\mathfrak{S}^{(1)}_6(C) = & p    \sum_{   b^2 C \preccurlyeq p^{3/2 + \alpha}  }  \sum_{   p \nmid \delta_0 \delta_1\delta_2  }  \sum_{\substack{\alpha \, (\text{mod } \delta_0)\\ \m_2\delta_2 + \epsilon_12^i b^2\bar{p} \delta_1 \alpha \equiv 0 \, (\text{mod } \delta_0)}}\sum_c \frac{W(c/C)}{ bc^{1/2}\delta_0\delta_1\delta_2}  \\
&
\times S\Big( \epsilon_2 2^j c^2\bar{p}, \frac{\m_2\delta_2 + \epsilon_12^i b^2\bar{p} \delta_1 \alpha}{\delta_0}, \delta_1\Big) S\Big(\m_1, \frac{\m_2\delta_2 \bar{\alpha} + \epsilon_12^i b^2\bar{p} \delta_1}{\delta_0}, \delta_2\Big) \\
& \times \Phi_{w_6}\Big(- \frac{\epsilon_2 2^jc^2\m_2\delta_2}{p\delta_0\delta_1^2}, - \frac{\epsilon_1 2^i b^2 \m_1\delta_1}{p\delta_0\delta_2^2} \Big). 
\end{split}
\end{displaymath}
For orientation, we should think of $\alpha$ as very small, and then the generic range is $C$ a bit larger than $p^{3/2}$ and $\delta_1$ a bit larger than $p$, all other variables essentially fixed. For this scenario, Poisson summation is beneficial, because it reduces the length of the $c$-sum. 
Note that the  $\alpha$-sum contains $\ll (\delta_0, \delta_1b^2, \m_2\delta_2)$ terms. 

We   apply Poisson summation to the $c$-sum and recast $\mathfrak{S}^{(1)}_6(C)$ as   \begin{displaymath}
\begin{split}
p    \sum_{   b^2 C \preccurlyeq p^{3/2 + \alpha}  } & \sum_{   p \nmid \delta_0 \delta_1\delta_2  }  \sum_{\substack{\alpha \, (\text{mod } \delta_0)\\ \m_2\delta_2 + \epsilon_12^i b^2\bar{p} \delta_1 \alpha \equiv 0 \, (\text{mod } \delta_0)}}  \sum_{\gamma \, (\text{mod } \delta_1)}  \sum_c e\Big( \frac{c\gamma}{\delta_1}\Big) \int \frac{W(x/C)}{ bx^{1/2}\delta_0\delta_1^2\delta_2} e\Big(- \frac{cx}{\delta_1}\Big)  \\
&
\times S\Big( \epsilon_2 2^j \gamma^2\bar{p}, \frac{\m_2\delta_2 + \epsilon_12^i b^2\bar{p} \delta_1 \alpha}{\delta_0}, \delta_1\Big) S\Big(\m_1, \frac{\m_2\delta_2 \bar{\alpha} + \epsilon_12^i b^2\bar{p} \delta_1}{\delta_0}, \delta_2\Big) \\
& \times  \Phi_{w_6}\Big(- \frac{\epsilon_2 2^jx^2\m_2\delta_2}{p\delta_0\delta_1^2}, - \frac{\epsilon_1 2^i b^2 \m_1\delta_1}{p\delta_0\delta_2^2} \Big) dx. 
\end{split}
\end{displaymath}
Repeated integration by parts together with \eqref{lem8b}   shows that the $x$-integral is negligible unless
$$c \preccurlyeq C^{\ast} = \frac{\delta_1}{C}\Big(1 + \frac{C(\m_2\delta_2)^{1/2}}{p^{1/2} \delta_0^{1/2} \delta_1} + \frac{C^{2/3} (b\m_2)^{1/3}\m_1^{1/6}}{\delta_1^{1/2} \delta_0^{1/2} p^{1/2}}\Big).$$
Next we consider the $\gamma$-sum: with $A = (\m_2\delta_2 + \epsilon_12^i b^2\bar{p} \delta_1 \alpha)/\delta_0$ we have 
\begin{displaymath}
\begin{split}
\sum_{\gamma \, (\text{mod } \delta_1)} & e\Big( \frac{c\gamma}{\delta_1}\Big)S ( \epsilon_2 2^j \gamma^2\bar{p}, A, \delta_1 )   = \sum_{\gamma \, (\text{mod } \delta_1)}  \sum_{\substack{f\, (\text{mod } \delta_1)\\ (f, \delta_1) = 1}}  e\Big( \frac{c\gamma  + \bar{f} \epsilon_2 2^j \gamma^2\bar{p} + f A }{\delta_1}\Big)\\
& =  \sum_{\gamma \, (\text{mod } \delta_1)}  \sum_{\substack{f\, (\text{mod } \delta_1)\\ (f, \delta_1) = 1}}  e\Big( \frac{cf\gamma  + f\epsilon_2 2^j \gamma^2\bar{p} + f A }{\delta_1}\Big) =  \sum_{\gamma \, (\text{mod } \delta_1)} r_{\delta_1} (c\gamma p+ \epsilon_2 2^j \gamma^2   + Ap)\\
& \ll  \sum_{\gamma \, (\text{mod } \delta_1)}  (\delta_1, c\gamma p + \epsilon_2 2^j \gamma^2  + Ap) \ll \delta_1^{1 + \varepsilon}.
\end{split}\end{displaymath}
As before we use Weil's bound for the Kloosterman sum modulo $\delta_2$. We estimate the $x$-integral trivially, using the bound \eqref{lem8c}. Thus we are left with bounding 
 \begin{displaymath}
\begin{split}
& p    \sum_{   b^2 C \preccurlyeq p^{3/2 + \alpha}  }  \sum_{ \substack{  \delta_0\delta_1  \preccurlyeq (C^4b^2\m_2^2\m_1)^{1/3}  p^{-1}\\ \delta_0 \delta_2 \ll (b^4\m_1^2C^2\m_2)^{1/3}p^{-1} 
}}  (\delta_0, \delta_1b^2,\m_2 \delta_2)   \sum_{c\ll C^{\ast}} C^{1/2}  \frac{ (\delta_2, \m_1)^{1/2}}{ b\delta_0\delta_1\delta_2^{1/2} }\Big(\frac{C^2   b^2 \m_1\m_2  }{p^2 \delta_0^2 \delta_1 \delta_2 }\Big)^{-1/4}\\
\preccurlyeq &p^{3/2} \Big(\frac{\m_1}{\m_2}\Big)^{1/4} \sum_{   b^2 C \preccurlyeq p^{3/2 + \alpha}  }  \sum_{ \substack{ \delta_0\delta_1  \preccurlyeq (C^4b^2\m_2^2\m_1)^{1/3}  p^{-1}\\ \delta_0 \delta_2\preccurlyeq (b^4\m_1^2C^2\m_2)^{1/3}p^{-1}  
}} \\
& \quad\quad\quad \Big(1 + \frac{\delta_1}{C} + \frac{(\m_2\delta_2)^{1/2}}{p^{1/2} \delta_0^{1/2} } + \frac{\delta_1^{1/2}(b\m_2)^{1/3}\m_1^{1/6}}{C^{1/3}  \delta_0^{1/2} p^{1/2}}\Big)  \frac{  (\delta_0, \delta_1b^2, \m_2\delta_2) }{ b^{3/2}\delta_0^{1/2}\delta_1^{3/4}\delta_2^{1/4} }.
\end{split}
\end{displaymath}
Summing over $\delta_1, \delta_2, \delta_0$, we obtain
\begin{displaymath}
\begin{split}
 & p^{3/2}\Big(\frac{\m_1}{\m_2}\Big)^{1/4}  \sum_{   b^2 C \preccurlyeq p^{3/2 + \alpha}  }  \frac{1}{b^{3/2}}\Big[\Big( \frac{(C^4b^2\m_2^2\m_1)^{1/3} }{ p }   \Big)^{1/4}\Big(\frac{(b^4\m_1^2C^2\m_2)^{1/3}}{p}\Big)^{3/4} \\
 &+ \frac{1}{C} \Big( \frac{(C^4b^2\m_2^2\m_1)^{1/3} }{ p }   \Big)^{5/4}\Big(\frac{(b^4\m_1^2C^2\m_2)^{1/3}}{p}\Big)^{3/4}  + \frac{\m_2^{1/2}}{p^{1/2}} \Big( \frac{(C^4b^2\m_2^2\m_1)^{1/3} }{ p }   \Big)^{1/4}\Big(\frac{(b^4\m_1^2C^2\m_2)^{1/3}}{p}\Big)^{5/4} \\
 & +  \frac{ (b\m_2)^{1/3}\m_1^{1/6}}{C^{1/3}   p^{1/2}}\Big(\frac{(C^4b^2\m_2^2\m_1)^{1/3} }{ p }   \Big)^{3/4}\Big(\frac{(b^4\m_1^2C^2\m_2)^{1/3}}{p}\Big)^{3/4}  \Big].
\end{split}
\end{displaymath}
Simplifying, we obtain
\begin{equation}\label{S6}
\begin{split}
 \mathfrak{S}^{(1)}_6(C) & \preccurlyeq p^{3/2}\Big(\frac{\m_1}{\m_2}\Big)^{1/4}  \sum_{   b^2 C \preccurlyeq p^{3/2 + \alpha}  }  \Big(\frac{C^{5/6} \m_1^{7/12} \m_2^{5/12}}{b^{1/3}p} + \frac{C^{7/6} \m_1^{11/12} \m_2^{13/12} b^{1/3}}{p^2} \Big)\\
 & \preccurlyeq  p^{3/2}\Big(\frac{\m_1}{\m_2}\Big)^{1/4} \Big(\frac{ (p^{3/2 + \alpha})^{5/6} \m_1^{7/12} \m_2^{5/12}}{ p} + \frac{ (p^{3/2 + \alpha})^{7/6} \m_1^{11/12} \m_2^{13/12}  }{p^2} \Big)\\
 & \preccurlyeq  \big(p^{7/4 + 5\alpha/6} \m_1^{5/6} \m_2^{1/6} +  p^{5/4 + 7\alpha/6}\big)\m_1^{7/6} \m_2^{5/6}.
 \end{split}
\end{equation}
 Collecting the estimates \eqref{D}, \eqref{s4}, \eqref{S5} and \eqref{S6}, choosing $\alpha  < 1/20$ and substituting this into the right hand side of \eqref{formula} we obtain the final asymptotic formula
 \begin{displaymath}
 c_{\text{cusp}} \sum_{\pi \in \mathcal{B}(p)} \frac{A_{\pi}(\m_1, \m_2) L^{\ast}(1/2, \pi, \text{sym}^2) h(\mu_{\pi})}{L(1, \pi, \text{Ad})} = \delta_{\substack{\m_1 = \square\\ \m_2 = \square}} \frac{Cp^2}{(\m_1\m_2^2)^{1/4}} + O((\m_1\m_2)^{5/2} p^{9/5}).
 \end{displaymath}
 provided that \eqref{m}  holds. 
This completes the proof of Theorem \ref{thm1}.

 \section{Proof  of the Corollary \ref{cor3}}

 
 We follow the method of \cite{RS}. Let $h$ be as in \eqref{testh},  $k \in \Bbb{N}$  and 
  $x = p^{\eta}$ for some very small constant $0 < \eta < 1/k$. With
 $$\mathcal{A}(\pi) = \sum_{2 \nmid m \leq x} \frac{A_{\pi}(m^2, 1)}{m^{1/2}}$$
 we  define
 \begin{displaymath}
 \begin{split}
  & S_1 := \sum_{\pi \in  \mathcal{B}(p)}  \frac{L^{\ast}(1/2, \pi, \text{{\rm sym}}^2)}{L(1, \pi, \text{{\rm Ad}})} h(\mu_{\pi})  \mathcal{A}(\pi)^{k} \overline{\mathcal{A}(\pi)}^{k-1}, \quad S_2 := \sum_{\pi \in  \mathcal{B}(p)} \frac{|\mathcal{A}(\pi)|^{2k}  h(\mu_{\pi}) }{L(1, \pi, \text{{\rm Ad}})}.
 \end{split}
 \end{displaymath} 
 In the following we regard $k$ as fixed and all implied constants may depend on $k$. 
 
 \begin{lemma}\label{s2upper}
 We have $S_2 \ll p^2 (\log p)^{k^2}$. 
 \end{lemma}

\textbf{Proof.}  
We compute $\mathcal{A}(\pi)^k$ by repeated application of \eqref{hecke1}, which shows
 \begin{equation}\label{prod}
\begin{split}
\prod_{j=1}^k A_{\pi}(m_j^2, 1) = &\sum_{\substack{\alpha_{2, 0} \alpha_{2, 1} \alpha_{2, 2} = m_2^2\\ \alpha_{2, 1} \mid m_1^2 \\ \alpha_{2, 2} \mid 1}} \cdots \sum_{\substack{\alpha_{j, 0}\alpha_{j, 1}\alpha_{j, 2}= m_j^2\\ \alpha_{j, 1} \mid \frac{m_1^2 \cdots m_{j-1}^2}{\alpha_{2, 1}^2 \alpha_{2, 2} \cdots \alpha_{j-1, 1}^2 \alpha_{j-1, 2}}\\ \alpha_{j, 2} \mid \frac{\alpha_{2, 1}\cdots \alpha_{j-1, 1}}{\alpha_{2, 2} \cdots \alpha_{j-1, 2}}}} \cdots \sum_{\substack{\alpha_{k, 0}\alpha_{k, 1}\alpha_{k, 2}= m_k^2\\ \alpha_{k, 1} \mid \frac{m_1^2 \cdots m_{k-1}^2}{\alpha_{2, 1}^2 \alpha_{2, 2} \cdots \alpha_{k-1, 1}^2 \alpha_{k-1, 2}}\\ \alpha_{k, 2} \mid \frac{\alpha_{2, 1}\cdots \alpha_{k-1, 1}}{\alpha_{2, 2} \cdots \alpha_{k-1, 2}}}} \\
&A_{\pi}\Big( \frac{m_1^2 \cdots m_k^2}{\alpha_{2, 1}^2 \alpha_{2, 2} \cdots \alpha_{k, 1}^2\alpha_{k, 2}}, \frac{\alpha_{2, 1}\cdots \alpha_{k, 1}}{\alpha_{2, 2} \cdots \alpha_{k, 2}}\Big)
\end{split}
\end{equation}
 for $(m_j, p) = 1$. We multiply this expression by its conjugate and sum over odd $m_j \leq x$ to obtain $|\mathcal{A}(\pi)|^{2k}$. We then apply the Kuznetsov formula to $S_2$, and as in the proof of Theorem \ref{thm1} we need to add the rest of spectrum in the sum over $\pi \in \mathcal{B}(p)$. Here we can do this simply by positivity since cuspidal oldforms and Eisenstein series $\pi$ satisfy the same (unramified) Hecke relations, so the corresponding summands are $|\mathcal{A}(\pi)|^{2k}  h(\mu_{\pi})/L(1, \pi, \text{Ad})$ and hence non-negative. If $\eta$ is sufficiently small, it is easy to see from \eqref{lem8a} and \eqref{lem8b} with $A$ sufficiently large that all off-diagonal terms in the Kuznetsov formula are negligible. Thus
\begin{displaymath}
\begin{split}
S_2 \ll &p^2\sum_{2 \nmid a_1, \ldots, a_k \leq x} \sum_{\substack{\alpha_{2, 0} \alpha_{2, 1} \alpha_{2, 2} = a_2^2\\ \alpha_{2, 1} \mid a_1^2 \\ \alpha_{2, 2} \mid 1}} \cdots \sum_{\substack{\alpha_{j, 0}\alpha_{j, 1}\alpha_{j, 2}= a_j^2\\ \alpha_{j, 1} \mid \frac{a_1^2 \cdots a_{j-1}^2}{\alpha_{2, 1}^2 \alpha_{2, 2} \cdots \alpha_{j-1, 1}^2 \alpha_{j-1, 2}}\\ \alpha_{j, 2} \mid \frac{\alpha_{2, 1}\cdots \alpha_{j-1, 1}}{\alpha_{2, 2} \cdots \alpha_{j-1, 2}}}} \cdots \sum_{\substack{\alpha_{k, 0}\alpha_{k, 1}\alpha_{k, 2}= a_k^2\\ \alpha_{k, 1} \mid \frac{a_1^2 \cdots a_{k-1}^2}{\alpha_{2, 1}^2 \alpha_{2, 2} \cdots \alpha_{k-1, 1}^2 \alpha_{k-1, 2}}\\ \alpha_{k, 2} \mid \frac{\alpha_{2, 1}\cdots \alpha_{k-1, 1}}{\alpha_{2, 2} \cdots \alpha_{k-1, 2}}}}\\
&  \sum_{2 \nmid b_1, \ldots, b_k \leq x} \sum_{\substack{\beta_{2, 0} \beta_{2, 1} \beta_{2, 2} = b_2^2\\ \beta_{2, 1} \mid b_1^2 \\ \beta_{2, 2} \mid 1}} \cdots \sum_{\substack{\beta_{j, 0}\beta_{j, 1}\beta_{j, 2}= b_j^2\\ \beta_{j, 1} \mid \frac{b_1^2 \cdots b_{j-1}^2}{\beta_{2, 1}^2 \beta_{2, 2} \cdots \beta_{j-1, 1}^2 \beta_{j-1, 2}}\\ \beta_{j, 2} \mid \frac{\beta_{2, 1}\cdots \beta_{j-1, 1}}{\beta_{2, 2} \cdots \beta_{j-1, 2}}}} \cdots \sum_{\substack{\beta_{k, 0}\beta_{k, 1}\beta_{k, 2}= b_k^2\\ \beta_{k, 1} \mid \frac{b_1^2 \cdots b_{k-1}^2}{\beta_{2, 1}^2 \beta_{2, 2} \cdots \beta_{k-1, 1}^2 \beta_{k-1, 2}}\\ \beta_{k, 2} \mid \frac{\beta_{2, 1}\cdots \beta_{k-1, 1}}{\beta_{2, 2} \cdots \beta_{k-1, 2}}}}\\
&\delta_{\frac{a_1^2 \cdots a_k^2}{\alpha_{2, 1}^2 \alpha_{2, 2} \cdots \alpha_{k, 1}^2\alpha_{k, 2}} = \frac{b_1^2 \cdots b_k^2}{\beta_{2, 1}^2 \beta_{2, 2} \cdots \beta_{k, 1}^2\beta_{k, 2}}} \delta_{\frac{\alpha_{2, 1}\cdots \alpha_{k, 1}}{\alpha_{2, 2} \cdots \alpha_{k, 2}} = \frac{\beta_{2, 1}\cdots \beta_{k, 1}}{\beta_{2, 2} \cdots \beta_{k, 2}}}
\frac{1}{\sqrt{a_1 \cdots a_k b_1 \cdots b_k}}.
\end{split}
\end{displaymath}
By positivity, we drop the condition that $a_j, b_j$ are odd, and we enlarge the summation condition $a_j, b_j \leq x$ to $a_1\cdots a_kb_1 \cdots b_k \leq x^{2k}$ which we encode by a smooth weight function $W(a_1\cdots a_kb_1 \cdots b_k/x^{2k})$ where $W = 1$ on $[0, 1]$ and $W = 0$ on $[2, \infty)$. By Mellin inversion, we obtain
\begin{equation}\label{s2}
S_2 \ll p^2 \int_{(2)} \widehat{W}(s) D(s) x^{2ks} \frac{ds}{2\pi i}
\end{equation}
where $\widehat{W}$ has a simple pole at 0 and 
\begin{displaymath}
\begin{split}
D(s) &=   \sum_{  a_1, \ldots, a_k } \sum_{\substack{\alpha_{2, 0} \alpha_{2, 1} \alpha_{2, 2} = a_2^2\\ \alpha_{2, 1} \mid a_1^2 \\ \alpha_{2, 2} \mid 1}} \cdots \sum_{\substack{\alpha_{j, 0}\alpha_{j, 1}\alpha_{j, 2}= a_j^2\\ \alpha_{j, 1} \mid \frac{a_1^2 \cdots a_{j-1}^2}{\alpha_{2, 1}^2 \alpha_{2, 2} \cdots \alpha_{j-1, 1}^2 \alpha_{j-1, 2}}\\ \alpha_{j, 2} \mid \frac{\alpha_{2, 1}\cdots \alpha_{j-1, 1}}{\alpha_{2, 2} \cdots \alpha_{j-1, 2}}}} \cdots \sum_{\substack{\alpha_{k, 0}\alpha_{k, 1}\alpha_{k, 2}= a_k^2\\ \alpha_{k, 1} \mid \frac{a_1^2 \cdots a_{k-1}^2}{\alpha_{2, 1}^2 \alpha_{2, 2} \cdots \alpha_{k-1, 1}^2 \alpha_{k-1, 2}}\\ \alpha_{k, 2} \mid \frac{\alpha_{2, 1}\cdots \alpha_{k-1, 1}}{\alpha_{2, 2} \cdots \alpha_{k-1, 2}}}}\\
&  \sum_{  b_1, \ldots, b_k  } \sum_{\substack{\beta_{2, 0} \beta_{2, 1} \beta_{2, 2} = b_2^2\\ \beta_{2, 1} \mid b_1^2 \\ \beta_{2, 2} \mid 1}} \cdots \sum_{\substack{\beta_{j, 0}\beta_{j, 1}\beta_{j, 2}= b_j^2\\ \beta_{j, 1} \mid \frac{b_1^2 \cdots b_{j-1}^2}{\beta_{2, 1}^2 \beta_{2, 2} \cdots \beta_{j-1, 1}^2 \beta_{j-1, 2}}\\ \beta_{j, 2} \mid \frac{\beta_{2, 1}\cdots \beta_{j-1, 1}}{\beta_{2, 2} \cdots \beta_{j-1, 2}}}} \cdots \sum_{\substack{\beta_{k, 0}\beta_{k, 1}\beta_{k, 2}= b_k^2\\ \beta_{k, 1} \mid \frac{b_1^2 \cdots b_{k-1}^2}{\beta_{2, 1}^2 \beta_{2, 2} \cdots \beta_{k-1, 1}^2 \beta_{k-1, 2}}\\ \beta_{k, 2} \mid \frac{\beta_{2, 1}\cdots \beta_{k-1, 1}}{\beta_{2, 2} \cdots \beta_{k-1, 2}}}}\\
&\delta_{\frac{a_1^2 \cdots a_k^2}{\alpha_{2, 1}^2 \alpha_{2, 2} \cdots \alpha_{k, 1}^2\alpha_{k, 2}} = \frac{b_1^2 \cdots b_k^2}{\beta_{2, 1}^2 \beta_{2, 2} \cdots \beta_{k, 1}^2\beta_{k, 2}}} \delta_{\frac{\alpha_{2, 1}\cdots \alpha_{k, 1}}{\alpha_{2, 2} \cdots \alpha_{k, 2}} = \frac{\beta_{2, 1}\cdots \beta_{k, 1}}{\beta_{2, 2} \cdots \beta_{k, 2}}}
\frac{1}{ (a_1 \cdots a_k b_1 \cdots b_k)^{1/2 +  s}}.
\end{split}
\end{displaymath}
By multiplicativity, the Dirichlet series factorizes into an Euler product
$$D(s) = \prod_p \Big(1 + \frac{c_1}{p^{1/2+ s}} + \frac{c_2}{p^{2(1/2+ s)}}   + \ldots\Big).$$

We clearly have $c_1 = 0$: assume that $a_j = p$ for some $1 \leq j \leq k$, say, and $a_i = 1$ for $i \not = j$ and $b_i = 1$ for all $1 \leq i \leq k$.  Then $\beta_{i, \ell} = 1$ for all $i, \ell$ and $\alpha_{i, \ell} = 1$ for all $\ell$ and $i \not= j$. Thus the second $\delta$-symbol implies $\alpha_{j, 1} = \alpha_{j, 2}$, but this contradicts the first $\delta$-symbol,

Next we claim that $c_2 = k^2$. Indeed, arguing as before, we see that the $a_j = p^2$ for some $j$, $a_i = 1$ for $i \not = j$, $b_i = 1$ for all $i$ is impossible and similarly $a_i = a_j = p$ for some $i \not =j$, $a_{\ell} = 1$ for $\ell \not\in \{i, j\}$, $b_\ell =1$ for all $\ell$ is impossible. Hence the only option is $a_i = p$, $b_j = p$ for one of the $k^2$ pairs $(i, j)$ and all other $a$'s and $b$'s are 1. The divisibility conditions imply $\alpha_{i, 1} = \alpha_{i, 2} = \beta_{j, 1} = \beta_{j, 2} = 1$, and hence $\alpha_{i, 0} = \beta_{j, 0} = p$, in other words all $\alpha$'s and $\beta$'s are determined.
  
We trivially have the coarse bound $c_r \leq (r+1)^{6k}$, since the $6k$ variables $a_i,  \alpha_{i, 1}, \alpha_{i, 2}, b_i,   \beta_{i, 1}, \beta_{i, 2}$ can take at most the $r+1$ values $1, p, p^2, \ldots, p^{r}$ and then $\alpha_{i, 0}, \beta_{i, 0}$ are determined.  We conclude
$$D(s) = \prod_{p} \Big(1 + \frac{k^2}{p^{1 + 2s}} + O\Big( \sum_{r=3}^{\infty} \frac{(r+1)^{6k}}{p^{r(1/2 + \Re s)}}\Big)\Big) = \zeta(1+2s)^{k^2} H(s)$$
where $H(s)$ is holomorphic in $\Re s > -1/6$.   Shifting the contour in \eqref{s2} past the pole of order $k^2 + 1$ at $s=0$ and observing that $\log x \asymp \log p$, we obtain the desired bound. 
  
  \begin{lemma}\label{s1lower}
 We have $S_1 \gg p^2 (\log p)^{k^2}$. 
 \end{lemma}    
     
\textbf{Proof.} This is similar as the  preceding proof, based on Theorem \ref{thm1}. We start from \eqref{prod} and repeatedly use \eqref{hecke1} with exchanged indices to arrive at
 \begin{displaymath}
\begin{split}
&\prod_{j=1}^k A_{\pi}(m_j^2, 1)\prod_{j=1}^{k-1} A_{\pi}(1, n_j^2) \\
= &\sum_{\substack{\alpha_{2, 0} \alpha_{2, 1} \alpha_{2, 2} = m_2^2\\ \alpha_{2, 1} \mid m_1^2 \\ \alpha_{2, 2} \mid 1}} \cdots \sum_{\substack{\alpha_{j, 0}\alpha_{j, 1}\alpha_{j, 2}= m_j^2\\ \alpha_{j, 1} \mid \frac{m_1^2 \cdots m_{j-1}^2}{\alpha_{2, 1}^2 \alpha_{2, 2} \cdots \alpha_{j-1, 1}^2 \alpha_{j-1, 2}}\\ \alpha_{j, 2} \mid \frac{\alpha_{2, 1}\cdots \alpha_{j-1, 1}}{\alpha_{2, 2} \cdots \alpha_{j-1, 2}}}} \cdots \sum_{\substack{\alpha_{k, 0}\alpha_{k, 1}\alpha_{k, 2}= m_k^2\\ \alpha_{k, 1} \mid \frac{m_1^2 \cdots m_{k-1}^2}{\alpha_{2, 1}^2 \alpha_{2, 2} \cdots \alpha_{k-1, 1}^2 \alpha_{k-1, 2}}\\ \alpha_{k, 2} \mid \frac{\alpha_{2, 1}\cdots \alpha_{k-1, 1}}{\alpha_{2, 2} \cdots \alpha_{k-1, 2}}}} \\
& \sum_{\substack{\beta_{1, 0}\beta_{1, 1}\beta_{1, 2} = n_1^2\\ \beta_{1, 1} \mid  \frac{m_1^2 \cdots m_k^2}{\alpha_{2, 1}^2 \alpha_{2, 2} \cdots \alpha_{k, 1}^2\alpha_{k, 2}}\\ \beta_{1, 2} \mid \frac{\alpha_{2, 1}\cdots \alpha_{k, 1}}{\alpha_{2, 2} \cdots \alpha_{k, 2}}}}  \cdots \sum_{\substack{\beta_{j, 0}\beta_{j, 1}\beta_{j, 2} = n_j^2\\ \beta_{j, 1} \mid  \frac{m_1^2 \cdots m_k^2}{\alpha_{2, 1}^2 \alpha_{2, 2} \cdots \alpha_{k, 1}^2\alpha_{k, 2}} \frac{\beta_{1, 2} \cdots \beta_{j-1, 2}}{\beta_{1, 1} \cdots \beta_{j-1, 1}}\\ \beta_{j, 2} \mid \frac{\alpha_{2, 1}\cdots \alpha_{k, 1}}{\alpha_{2, 2} \cdots \alpha_{k, 2}}  \frac{n_1^2 \cdots n_{j-1}^2}{\beta_{1, 1}\beta_{1, 2}^2 \cdots \beta_{j-1, 1}\beta_{j-1, 2}^2}}}  \cdots \sum_{\substack{\beta_{k-1, 0}\beta_{k-1, 1}\beta_{k-1, 2} = n_{k-1}^2\\ \beta_{k-1, 1} \mid  \frac{m_1^2 \cdots m_k^2}{\alpha_{2, 1}^2 \alpha_{2, 2} \cdots \alpha_{k, 1}^2\alpha_{k, 2}} \frac{\beta_{1, 2} \cdots \beta_{k-2, 2}}{\beta_{1, 1} \cdots \beta_{k-2, 1}}\\ \beta_{k-1, 2} \mid \frac{\alpha_{2, 1}\cdots \alpha_{k, 1}}{\alpha_{2, 2} \cdots \alpha_{k, 2}}  \frac{n_1^2 \cdots n_{k-2}^2}{\beta_{1, 1}\beta_{1, 2}^2 \cdots \beta_{k-2, 1}\beta_{k-2, 2}^2}}}\\
&A_{\pi}\Big( \frac{m_1^2 \cdots m_k^2}{\alpha_{2, 1}^2 \alpha_{2, 2} \cdots \alpha_{k, 1}^2\alpha_{k, 2}} \frac{\beta_{1, 2} \cdots \beta_{k-1, 2}}{\beta_{1, 1} \cdots \beta_{k-1, 1}}, \frac{\alpha_{2, 1}\cdots \alpha_{k, 1}}{\alpha_{2, 2} \cdots \alpha_{k, 2}}\frac{n_1^2 \cdots n_{k-1}^2}{\beta_{1, 1}\beta_{1, 2}^2 \cdots \beta_{k-1, 1}\beta_{k-1, 2}^2}\Big).
\end{split}
\end{displaymath}
We are now in a position to apply Theorem \ref{thm1}.  This gives us
\begin{displaymath}
\begin{split}
S_1 &= Cp^2 \sum_{2 \nmid m_1, \ldots, m_k \leq x} \sum_{2 \nmid n_1, \ldots, n_{k-1} \leq x} \frac{1}{(m_1\cdots m_k n_1 \cdots n_{k-1})^{1/2}}\\
&\sum_{\substack{\alpha_{2, 0} \alpha_{2, 1} \alpha_{2, 2} = m_2^2\\ \alpha_{2, 1} \mid m_1^2 \\ \alpha_{2, 2} \mid 1}} \cdots \sum_{\substack{\alpha_{j, 0}\alpha_{j, 1}\alpha_{j, 2}= m_j^2\\ \alpha_{j, 1} \mid \frac{m_1^2 \cdots m_{j-1}^2}{\alpha_{2, 1}^2 \alpha_{2, 2} \cdots \alpha_{j-1, 1}^2 \alpha_{j-1, 2}}\\ \alpha_{j, 2} \mid \frac{\alpha_{2, 1}\cdots \alpha_{j-1, 1}}{\alpha_{2, 2} \cdots \alpha_{j-1, 2}}}} \cdots \sum_{\substack{\alpha_{k, 0}\alpha_{k, 1}\alpha_{k, 2}= m_k^2\\ \alpha_{k, 1} \mid \frac{m_1^2 \cdots m_{k-1}^2}{\alpha_{2, 1}^2 \alpha_{2, 2} \cdots \alpha_{k-1, 1}^2 \alpha_{k-1, 2}}\\ \alpha_{k, 2} \mid \frac{\alpha_{2, 1}\cdots \alpha_{k-1, 1}}{\alpha_{2, 2} \cdots \alpha_{k-1, 2}}}} \\
& \sum_{\substack{\beta_{1, 0}\beta_{1, 1}\beta_{1, 2} = n_1^2\\ \beta_{1, 1} \mid  \frac{m_1^2 \cdots m_k^2}{\alpha_{2, 1}^2 \alpha_{2, 2} \cdots \alpha_{k, 1}^2\alpha_{k, 2}}\\ \beta_{1, 2} \mid \frac{\alpha_{2, 1}\cdots \alpha_{k, 1}}{\alpha_{2, 2} \cdots \alpha_{k, 2}}}}  \cdots \sum_{\substack{\beta_{j, 0}\beta_{j, 1}\beta_{j, 2} = n_j^2\\ \beta_{j, 1} \mid  \frac{m_1^2 \cdots m_k^2}{\alpha_{2, 1}^2 \alpha_{2, 2} \cdots \alpha_{k, 1}^2\alpha_{k, 2}} \frac{\beta_{1, 2} \cdots \beta_{j-1, 2}}{\beta_{1, 1} \cdots \beta_{j-1, 1}}\\ \beta_{j, 2} \mid \frac{\alpha_{2, 1}\cdots \alpha_{k, 1}}{\alpha_{2, 2} \cdots \alpha_{k, 2}}  \frac{n_1^2 \cdots n_{j-1}^2}{\beta_{1, 1}\beta_{1, 2}^2 \cdots \beta_{j-1, 1}\beta_{j-1, 2}^2}}}  \cdots \sum_{\substack{\beta_{k-1, 0}\beta_{k-1, 1}\beta_{k-1, 2} = n_{k-1}^2\\ \beta_{k-1, 1} \mid  \frac{m_1^2 \cdots m_k^2}{\alpha_{2, 1}^2 \alpha_{2, 2} \cdots \alpha_{k, 1}^2\alpha_{k, 2}} \frac{\beta_{1, 2} \cdots \beta_{k-2, 2}}{\beta_{1, 1} \cdots \beta_{k-2, 1}}\\ \beta_{k-1, 2} \mid \frac{\alpha_{2, 1}\cdots \alpha_{k, 1}}{\alpha_{2, 2} \cdots \alpha_{k, 2}}  \frac{n_1^2 \cdots n_{k-2}^2}{\beta_{1, 1}\beta_{1, 2}^2 \cdots \beta_{k-2, 1}\beta_{k-2, 2}^2}}}\\
&\delta_{\substack{\alpha_{2, 2} \cdots \alpha_{k, 2} \beta_{1, 2} \cdots \beta_{k-1, 2}\beta_{1, 1} \cdots \beta_{k-1, 1} = \square\\ \alpha_{2, 1} \cdots \alpha_{k, 1} \alpha_{2, 2} \cdots \alpha_{k, 2}\beta_{1, 1} \cdots \beta_{k-1, 1} = \square}} \frac{(\alpha_{2, 2} \ldots \alpha_{k, 2}\beta_{1, 2} \ldots \beta_{k-1, 2}\beta_{1, 1} \ldots \beta_{k-1, 1})^{3/4}}{(m_1 \cdots m_k)^{1/2} (n_1 \ldots n_{k-1})} + O(p^{2 - \delta} x^{O(1)})
\end{split}
\end{displaymath}
with $C$ as in Theorem \ref{thm1}.  By positivity, we can drop certain terms to obtain a lower bound. We restrict the sum to $$\alpha_{j, 1} = \alpha_{j, 2}  = 1\,\, (2 \leq j \leq k), \quad \beta_{j, 1} = n_{j}^2, \quad \beta_{j, 2} = 1\,\, (1 \leq j \leq k-1).$$
In this way we see
\begin{displaymath}
\begin{split}
S_1 &\gg  p^2 \underset{n_1^2 \cdots n_{k-1}^2 \mid m_1^2 \cdots m_k^2}{\sum_{2 \nmid m_1, \ldots, m_k \leq x} \sum_{2 \nmid n_1, \ldots, n_{k-1} \leq x}} \frac{1}{m_1\cdots m_k } + O(p^{2 - \delta} x^{O(1)}).
\end{split}
\end{displaymath}
We write $n_1 \cdots n_{k-1} d = m_1\cdots m_k$ and restrict the summation condition to $m_1 \cdots m_k \leq x$, so that
\begin{displaymath}
\begin{split}
S_1 &\gg  p^2\sum_{2 \nmid m_1 \cdots m_k \leq x} \frac{\tau_{k}(m_1\cdots m_k)}{m_1\cdots m_k } + O(p^{2 - \delta} x^{O(1)}) = p^2\sum_{2 \nmid m\leq  x} \frac{\tau_{k}(m)^2}{m } + O(p^{2 - \delta} x^{O(1)}) 
\end{split}
\end{displaymath}
where $\tau_k$ is the $k$-fold divisor function. If $x$ is a sufficiently small power of $p$, the main term dominates the error term and is $\asymp p^2 (\log p)^{k^2}$, as desired. \\

The proof of Corollary \ref{cor3} follows now directly from Lemma \ref{s2upper}, Lemma \ref{s1lower} and H\"older's inequality:
$$\sum_{\pi \in  \mathcal{B}(p)} \frac{|L^{\ast}(1/2,  \pi, \text{\rm sym}^2)|^{2k}}{L(1, \pi, \text{{\rm Ad}})} h(\mu_{\pi})  \geq \frac{S_1^{2k}}{S_2^{2k-1}} \gg p^2 (\log p)^{k^2}.$$

\end{document}